# ZARISKI DENSITY OF MONODROMY GROUPS VIA A PICARD–LEFSCHETZ TYPE FORMULA†

JINXING XU

Abstract. For the universal family of cyclic covers of projective spaces branched along hyperplane arrangements in general position, we consider its monodromy group acting on an eigenspace of the middle cohomology of the fiber. We prove the monodromy group is Zariski dense in the corresponding linear group. As an application, we show the fundamental group of the moduli space of hyperplane arrangements is large. It can be viewed as a degenerate analogy of Carlson-Toledo's result about the monodromy groups of smooth hypersurfaces [Duke Math. J. **97**(3) (1999): 621-648]. The main ingredient in the proof is a Picard-Lefschetz type formula for a suitable degeneration of this family.

## 1. Introduction

Given a family $\pi : \mathcal{X} \to S$ of complex projective varieties which is of relative dimension $n$ and topologically a locally trivial fibration, write $X = \pi^{-1}(s)$ for a typical fiber and consider the $\pi_1(S, s)$-action on the middle cohomology group $H^n(X,\ \mathbb{R})$. We define the *monodromy group* $\Gamma \subset Aut(H^n(X,\ \mathbb{R}))$ of the family $\pi$ to be the image of $\pi_1(S, s)$. The determination of the monodromy group is often both an important step in many problems and a very interesting question of its own. While usually it is difficult to determine the monodromy group itself, the so called *algebraic monodromy group*, defined as the smallest algebraic subgroup of $Aut(H^n(X,\ \mathbb{R}))$ containing $\Gamma$, is relatively easy to handle and we denote it by $Mon$. In general, we use $G^0$ to denote the connected component of an algebraic group $G$ containing the identity, and we call $Mon^0$ the *connected algebraic monodromy group* of the family $\pi$.

This paper is concerned with the computation of $Mon^0$ for some naturally arising algebraic families. To state an interesting special case of our main result Theorem 2.1, suppose $m \geq n + 3$ is an even positive integer and let $\mathfrak{M}$ be the coarse moduli space of ordered $m$ hyperplane arrangements in $\mathbb{P}^n$ in general position (cf. subsection 2.1). Let $\mathcal{X} \xrightarrow{f} \mathfrak{M}$ be the universal family of double covers of $\mathbb{P}^n$ branched along $m$ hyperplanes in general position. For a typical fiber $X := f^{-1}(s)$, let $H^n(X,\ \mathbb{R})_{(1)}$ be the

†This version reformats the author manuscript dated December 24, 2016. Its mathematical content is reproduced without incorporating the correction below. A separate corrigendum correcting the proof of Proposition 4.2 is included with this submission. The statements of the main theorems are unchanged.

$-1$-eigenspace of $H^n(X,\ \mathbb{R})$ under the natural action of $\mathbb{Z}/2\mathbb{Z}$. We have the obvious decomposition $H^n(X,\ \mathbb{R}) = H^n(\mathbb{P}^n,\ \mathbb{R}) \oplus H^n(X,\ \mathbb{R})_{(1)}$ and $\Gamma$ acts trivially on $H^n(\mathbb{P}^n,\ \mathbb{R})$. Moreover, $H^n(\mathbb{P}^n,\ \mathbb{R}) = 0$ if $n$ is odd. The cup product provides a non-degenerate paring $Q$ on $H^n(X,\ \mathbb{R})_{(1)}$, symmetric or alternating according to the parity of $n$ (cf. subsection 2.3). We will prove:

**Theorem 1.1.** The connected algebraic monodromy group of the family $f$ is:

$$Mon^0 = Aut^0(H^n(X,\ \mathbb{R})_{(1)}, Q) = \begin{cases} Sp(H^n(X,\ \mathbb{R}), Q), & \text{if } n \text{ is odd;} \\ SO(H^n(X,\ \mathbb{R})_{(1)}, Q), & \text{if } n \text{ is even.} \end{cases}$$

In other words, $\Gamma$ is Zariski dense in $Aut^0(H^n(X,\ \mathbb{R})_{(1)}, Q)$.

If $n = 1$, we get the universal family of hyperelliptic curves, and in this case our Theorem 1.1 follows from a result of A'Campo [1], who proved that $\Gamma$ is in fact of finite index in the integral symplectic group $Sp(H^1(X, \mathbb{Z}), Q)$. Motivated by a conjecture of Dolgachev, Gerkmann-Sheng-Van Straten-Zuo [11] treated the $n = 3$, $m = 8$ case.

Although the fiber $X = f^{-1}(s)$ of the family $\mathcal{X} \xrightarrow{f} \mathfrak{M}$ is singular, we can show $X$ is a finite group quotient of a smooth projective variety (cf. subsection 2.2). Hence $H^n(X, \mathbb{Q})_{(1)}$ underlies a natural $\mathbb{Q}$-polarized Hodge structure of weight $n$, and we can define the *special Mumford-Tate group* $MT(X)$ as the smallest $\mathbb{Q}$-subgroup of $Aut(H^n(X, \mathbb{Q})_{(1)})$ whose real points containing the image of the circle group $U(1)$. By results in [5], for a very general $s \in \mathfrak{M}$, the connected algebraic monodromy group $Mon^0$ is contained in the derived group $MT(X)^{der}$ of $MT(X)$. Then the above theorem implies:

**Corollary 1.2.** For a very general $s \in \mathfrak{M}$, let $X = f^{-1}(s)$. Then the special Mumford-Tate group of $X$ is

$$MT(X) = Aut^0(H^n(X,\ \mathbb{Q})_{(1)}, Q) = \begin{cases} Sp(H^n(X,\ \mathbb{Q}), Q), & \text{if } n \text{ is odd;} \\ SO(H^n(X,\ \mathbb{Q})_{(1)}, Q), & \text{if } n \text{ is even.} \end{cases}$$

Concerning the monodromy group itself, it is natural to ask:

**Question 1.1.** For which pairs of positive integers $(n, m)$, the monodromy group $\Gamma$ is of finite index in the group $Aut(H^n(X,\ \mathbb{Z})_{(1)}, Q)$ ?

In view of the results of A'Campo in [1], we conjecture that the answer is yes for each even integer $m \geq n + 3$. See also [23, 24] for the arithmeticity results of monodromy groups arising from higher degree cyclic covers of $\mathbb{P}^1$.

**Strategy of the proof:** A standard technique to show the Zariski density of the monodromy group is Deligne's criterion (cf. Deligne [6], section 4.4). We also use this criterion to prove Theorem 1.1, but not in a straightforward way. Below we explain it in more detail. For ease of notations, we assume $n$ is odd. We will introduce a complex algebraic group $\widetilde{Mon}(\mathbb{C}) \subset Sp(H^n(X,\ \mathbb{C}), Q)$ in subsection 4.1 and reduce the problem to proving $\widetilde{Mon}(\mathbb{C}) = Sp(H^n(X,\ \mathbb{C}), Q)$.

The main input in Deligne's criterion is a spanning set $R$ of "vanishing cycles" in $H^n(X,\ \mathbb{C})$ such that $\widetilde{Mon}(\mathbb{C})$ is the smallest algebraic subgroup of $Sp(H^n(X,\ \mathbb{C}), Q)$ containing the transvections generated by $R$, and $R$ consists of a single $\widetilde{Mon}(\mathbb{C})$-orbit. In practice, the vanishing cycles can be constructed by Lefschetz degenerations and the condition that $R$ consists of a single $\widetilde{Mon}(\mathbb{C})$-orbit can be guaranteed by the irreducibility of the discriminant locus. In our situation, we can produce a spanning subset $R \subset H^n(X,\ \mathbb{C})$ by considering certain degenerations in an analogous way of the Lefschetz degenerations. But we can not show $R$ consists of a single $\widetilde{Mon}(\mathbb{C})$-orbit as usually done, since the discriminant locus of the parameter space is reducible.

In order to proceed, we slightly generalise Deligne's criterion in subsection 4.2 to deduce that, under the weaker assumption that $R$ consists of a finite number of $\widetilde{Mon}(\mathbb{C})$-orbits, there are only two possibilities:

(i) $\widetilde{Mon}(\mathbb{C}) = Sp(H^n(X,\ \mathbb{C}), Q)$;

(ii) there exists a nontrivial $\widetilde{Mon}(\mathbb{C})$-invariant linear subspace $U$ of $H^n(X,\ \mathbb{C})$ such that $R \subset U \cup U^{\perp}$.

To exclude case (ii), we will show there is a sequence of group embeddings

$$\widetilde{Mon}(\mathbb{C}) \hookrightarrow \overline{Mon} \hookrightarrow Sp(H^n(X,\ \mathbb{C}), Q),$$

satisfying:

(1) every $\widetilde{Mon}(\mathbb{C})$-invariant linear subspace of $H^n(X,\ \mathbb{C})$ is also $\overline{Mon}$-invariant;

(2) $\forall\ \alpha \neq 0 \in R$, the orbit $\overline{Mon} \cdot \alpha$ linearly spans $H^n(X,\ \mathbb{C})$.

So if a $\widetilde{Mon}(\mathbb{C})$-invariant subspace $U$ of $H^n(X,\ \mathbb{C})$ contains a nonzero element of $R$, then $U = H^n(X,\ \mathbb{C})$. A similar argument replacing $U$ by $U^{\perp}$ excludes the case (ii), and we are done.

**Organisation of the paper:** In Section 2, we introduce the family of cyclic covers of $\mathbb{P}^n$ branched along hyperplane arrangements in general position and state the main result Theorem 2.1. We relate these cyclic covers with smooth complete intersections in projective spaces, using the construction of Kummer covers. Then we relate certain cyclic covers of $\mathbb{P}^n$ branched along hyperplane arrangements with cyclic covers of $\mathbb{P}^1$ branched at some points, and in this way we obtain some Hodge theoretic information of $H^n(X,\ \mathbb{Q})$.

In Section 3, we derive a Picard-Lefschetz type formula for a one-parameter degeneration of Fermat type complete intersections, based on F. Pham's generalised Picard-Lefschetz formula. Using the construction of Kummer covers, we get a Picard-Lefschetz type formula for a one-parameter degeneration of cyclic covers of $\mathbb{P}^n$ branched along hyperplane arrangements.

In Section 4, we give the proofs of Theorem 1.1 and Theorem 2.1. As explained above, the main algebraic tool is a generalised version of Deligne's criterion, and we generalise both Deligne's original version and Carlson-Toledo's complex reflection version in [3]. In order to treat a degenerate case where we can not use these generalised versions of Deligne's criterion, we also include an alternative type result of semi-simple Lie groups from [20].

**Notation:** Throughout this paper, we use the following notations:

- We fix a generator $\sigma$ of the cyclic group $\mathbb{Z}/r\mathbb{Z}$ and a primitive $r$-th root of unit $\zeta_r$. Let $M$ be a $\mathbb{C}$-linear space, or a sheaf of $\mathbb{C}$-linear spaces on a scheme, on which the group $\mathbb{Z}/r\mathbb{Z} =< \sigma >$ acts. For $i \in \mathbb{Z}/r\mathbb{Z}$ we write $M_{(i)} := \{x \in M \mid \sigma(x) = \zeta_r^i x\}$, which in the sheaf case has to be interpreted on the level of local sections. We refer to $M_{(i)}$ as the $i$-eigenspace of $M$. Note $M = \oplus_{i\in\mathbb{Z}/r\mathbb{Z}} M_{(i)}$.
- Suppose $V$ is a real or complex linear space of dimension $n$, and $Q$ a nondegenerate bilinear form on $V$. We use $Aut(V, Q)$ to denote the (real or complex) algebraic group of linear isomorphisms of $V$ preserving $Q$, and $PAut(V, Q)$ to denote the adjoint group of $Aut(V, Q)$ (that is, $Aut(V, Q)$ modulo its center). If $Q$ is an alternating form, then $Aut(V, Q) = Sp(V, Q)$ is the symplectic group and $PAut(V, Q) = Sp(V, Q)/\{\pm 1\}$. If $Q$ is a symmetric form, then $Aut(V, Q) = O(V, Q)$ is the orthogonal group and $PAut(V, Q) = O(V, Q)/\{\pm 1\}$.
- Suppose $V$ is a complex linear space of finite dimension, and $H$ a nondegenerate Hermitian form on $V$. We use $Aut(V, H)$ or $U(V, H)$ to denote the algebraic group (defined over $\mathbb{R}$) of linear isomorphisms of $V$ preserving $H$, and $PAut(V, H)$ or $PU(V, H)$ to denote the adjoint group of $Aut(V, H)$ (that is, $Aut(V, H)$ modulo its center). Explicitly, $PAut(V, H) = U(V, H)/U(1)$.
- We use $\llcorner x \lrcorner$ to denote the greatest integer not exceeding $x$.

## 2. Geometric Constructions from Hyperplane Arrangements

Given a hyperplane arrangement $\mathfrak{A}$ in $\mathbb{P}^n$ in general position, the cyclic cover of $\mathbb{P}^n$ branched along $\mathfrak{A}$ is a natural geometric background to study many related problems, for example, the cohomology of a local system on the complement of $\mathfrak{A}$ (cf. [8]), the monodromy of Aomoto-Gelfand's generalised hypergeometric functions (cf. [2, 10]), etc. In this section, we recall the definition of these cyclic covers and study two constructions related to these cyclic covers.

2.1. **General set up and statement of the main result.** From now on to the end of this paper, we fix three positive integers $m, n, r$ such that $r$ divides $m$ and $m \geq n+3$. Given an ordered arrangement $\mathfrak{A} = (H_1, \cdots, H_m)$ of hyperplanes in $\mathbb{P}^n$, we say $\mathfrak{A}$ is in general position if no $n+1$ of the hyperplanes intersect in a point, or equivalently, if the divisor $\sum_{i=1}^m H_i$ has simple normal crossings. Denote the coarse moduli space of ordered $m$ hyperplane arrangements in $\mathbb{P}^n$ in general position by $\mathfrak{M}_{AR}$. It is an open subset of $(\mathbb{P}^n)^m/PGL(n+1)$. Since $r|m$, for each hyperplane arrangement $(H_1, \cdots, H_m)$ in $\mathfrak{M}_{AR}$, we can define a (unique up to isomorphism) degree $r$ cyclic cover of $\mathbb{P}^n$ branched along the divisor $\sum_{i=1}^m H_i$. In this way we obtain a universal family $\mathcal{X}_{AR} \xrightarrow{f} \mathfrak{M}_{AR}$ of degree $r$ cyclic covers of $\mathbb{P}^n$ branched along $m$ hyperplane arrangements in general position. Because the complements of any two hyperplane arrangements in general position are homeomorphic to each other, it follows this family is a locally trivial topological fibration, although it is not a smooth family (see [19] for an explicit crepant resolution algorithm). Hence for a typical fiber $X := f^{-1}(s)$, we have the monodromy group $\Gamma \subset Aut(H^n(X, \mathbb{R}))$ of $f$ acting on the middle cohomology of $X$.

As we will see from the construction in subsection 2.2, although $X$ is singular, the cup product still provides a non-degenerate bilinear form $Q$ on $H^n(X, \mathbb{R})$, symmetric or alternating according to the parity of $n$. Note that $Q$ is invariant under the action of $\Gamma$. The cyclic cover structure induces a natural action of the cyclic group $\mathbb{Z}/r\mathbb{Z}$ on $H^n(X, \mathbb{R})$. By our notation convention at the end of Section 1, we have the following ($\Gamma$-invariant) eigenspace decomposition:

$$H^n(X, \mathbb{C}) = \oplus_{i=0}^{r-1} H^n(X, \mathbb{C})_{(i)}.$$

It is easy to see $\Gamma$ acts trivially on $H^n(X, \mathbb{C})_{(0)} = H^n(\mathbb{P}^n, \mathbb{C})$. For $1 \leq i < \frac{r}{2}$, the subspace $H^n(X, \mathbb{C})_{(i)} \oplus H^n(X, \mathbb{C})_{(r-i)}$ underlies a natural real structure, i.e., there is a real subspace $V_{(i)} \subset H^n(X, \mathbb{R})$ such that $V_{(i)} \otimes \mathbb{C} = H^n(X, \mathbb{C})_{(i)} \oplus H^n(X, \mathbb{C})_{(r-i)}$. We remark in this case $V_{(i)} \simeq H^n(X, \mathbb{C})_{(i)}$ as $\mathbb{R}$-linear spaces. Similarly, if $r$ is an even integer, there is a real subspace $V_{(\frac{r}{2})} \subset H^n(X, \mathbb{R})$ such that $V_{(\frac{r}{2})} \otimes \mathbb{C} = H^n(X, \mathbb{C})_{(\frac{r}{2})}$. In any case, for each integer $i$ with $1 \leq i \leq \llcorner \frac{r}{2} \lrcorner$, let $Mon_{(i)}$ be the smallest real algebraic subgroup of $Aut(V_{(i)})$ containing the image of the monodromy representation

$$\rho_i : \pi_1(\mathfrak{M}_{AR}, s) \to Aut(V_{(i)}).$$

We call $Mon^0_{(i)}$ *the $i$-th eigen connected algebraic monodromy group* of the family $f$.

In subsection 2.3, we will show that the induced Hermitian form

$$H(\alpha, \beta) := \sqrt{-1}^n Q(\alpha, \bar{\beta})$$

on $H^n(X, \mathbb{C})$ is $\Gamma$-invariant, and for $1 \leq i \leq r-1$, its restriction on $H^n(X, \mathbb{C})_{(i)}$ is non-degenerate, with signature $(p_i, q_i)$ or $(q_i, p_i)$ :

$$p_i = \sum_{j=0}^{\llcorner \frac{n}{2} \lrcorner} \binom{m-ki-1}{n-2j} \binom{ki-1}{2j},$$

$$q_i = \sum_{j=0}^{\llcorner \frac{n-1}{2} \lrcorner} \binom{m-ki-1}{n-2j-1} \binom{ki-1}{2j+1}.$$

where $k = \frac{m}{r}$.

The main result of this paper is the following

**Theorem 2.1.** Given an integer $i$ satisfying $1 \leq i \leq \llcorner \frac{r}{2} \lrcorner$ and $mi \geq 2r$, the $i$-th eigen *connected algebraic monodromy group* $Mon^0_{(i)}$ of the universal family $\mathcal{X}_{AR} \xrightarrow{f} \mathfrak{M}_{AR}$ of degree $r$ cyclic covers of $\mathbb{P}^n$ branched along $m$ hyperplane ararngements in general position is:

(1) If $r = 2i$, then

$$Mon^0_{(\frac{r}{2})} = Aut^0(V_{(\frac{r}{2})}, Q) = \begin{cases} Sp(V_{(\frac{r}{2})}, Q), & \text{if } n \text{ is odd;} \\ SO(V_{(\frac{r}{2})}, Q), & \text{if } n \text{ is even.} \end{cases}$$

(2) If $1 \leq i < \frac{r}{2}$, then

$$Mon^0_{(i)} = SU(H^n(X,\ \mathbb{C})_{(i)}, H) \simeq SU(p_i, q_i).$$

From this theorem, we can deduce Theorem 1.1 by taking $r = 2$. The $n = 1$ case has been studied extensively (see e.g. [7, 14, 18]). In [20], the Calabi-Yau cases ($m = n + \frac{m}{r} + 1$ ) are proved, but the infinitesimal computations there rely heavily on the Calabi-Yau condition, so it can not be applied to the general case directly.

As an application of our main result, we can give a direct proof of the fact that the fundamental group of the coarse moduli space $\mathfrak{M}_{AR}$ is large. More precisely, we have

**Corollary 2.2.** The fundamental group $\pi_1(\mathfrak{M}_{AR}, s)$ is large, that is, there is a homomorphism of $\pi_1(\mathfrak{M}_{AR}, s)$ to a noncompact semisimple real algebraic group which has Zariski-dense image.

*Proof.* For each $m \geq n + 3$, it is elementary to see we can find appropriate $r$ and $i$ satisfying $r|m$, $1 \leq i \leq \llcorner \frac{r}{2} \lrcorner$ and $mi \geq 2r$. Then by definition, the monodromy representation

$$\rho_i : \pi_1(\mathfrak{M}_{AR}, s) \to Mon_{(i)}$$

is a homomorphism from $\pi_1(\mathfrak{M}_{AR}, s)$ to a real algebraic group with Zariski-dense image. Since a direct computation shows $p_i$, $q_i \geq 1$, we deduce from Theorem 2.1 that $Mon^0_{(i)}$ is noncompact semisimple and we are done.

□

**Remark 2.1.** This corollary first appeared in [20], where it is deduced using only double covers, and the homomorphisms from $\pi_1(\mathfrak{M}_{AR}, s)$ to semisimple real algebraic groups are not given explicitly. The proof above is more natural and more direct than the original one.

Note large groups always contain a free group of rank two by the Tits alternative [22]. In [3], Carlson and Toledo considered the universal family of degree $d$ smooth hypersurfaces in $\mathbb{P}^{n+1}$, and they showed the kernel of the monodromy representation is large if $d > 2$ and $(d, n) \neq (3, 1), (3, 0)$. Since hyperplane arrangements are degenerate hypersurfaces, and the monodromy acts on these degenerate hypersurfaces by permuting these hyperplane arrangements, we see the kernel of this monodromy representation is commensurable with $\pi_1(\mathfrak{M}_{AR}, s)$. So Corollary 2.2 can be viewed as a degenerate analogy of Carlson-Toledo's result.

2.2. **The Kummer cover.** Much of the material in this and the next subsection comes from [20]. Let $\mathfrak{A} = (H_1, \cdots, H_m) \in \mathfrak{M}_{AR}$ be an ordered arrangement in general position. It is easy to verify that under the automorphism group of $\mathbb{P}^n$ one can transform in a unique way the ordered first $n+2$ hyperplanes of $\mathfrak{A}$ into the ordered $n+2$ hyperplanes in $\mathbb{P}^n$, that are given by the first $n+2$ columns in the following $(n+1) \times m$ matrix:

$$B = (b_{ij}) := \begin{pmatrix} 1 & 0 & \cdots & 0 & 1 & 1 & \cdots & 1 \\ 0 & 1 & & 0 & 1 & a_{11} & \cdots & a_{1,m-n-2} \\ \vdots & & \ddots & \vdots & \vdots & \vdots & & \vdots \\ 0 & & & 1 & 1 & a_{n1} & \cdots & a_{n,m-n-2} \end{pmatrix}$$

Here the $j$-th column corresponds to the defining equation

$$\sum_{i=0}^{n} b_{ij} x_i = 0$$

of the hyperplane $H_j$, and $[x_0 : \cdots : x_n]$ are the homogeneous coordinates on $\mathbb{P}^n$.

Given $\mathfrak{A}$ as described above, let $Y$ be the complete intersection of the $m-n-1$ hypersurfaces in $\mathbb{P}^{m-1}$ defined by the equations:

$$y_{n+1}^r - (y_0^r + \sum_{j=1}^{n} y_j^r) = 0;$$

$$y_{n+i+1}^r - (y_0^r + \sum_{j=1}^{n} a_{ji} y_j^r) = 0, \ 1 \le i \le m-n-2.$$

Here $[y_0 : \cdots : y_{m-1}]$ are the homogeneous coordinates on $\mathbb{P}^{m-1}$. Since $\mathfrak{A}$ is in general position, the space $Y$ is smooth (see Proposition 3.1.2 in [21]). We call $Y$ the Kummer cover of the degree $r$ cyclic cover $X$ of $\mathbb{P}^n$ branched along $\mathfrak{A}$.

Let $N = \oplus_{j=0}^{m-1} \mathbb{Z}/r\mathbb{Z}$. Consider the following group

$$N_1 := Ker(N \to \mathbb{Z}/r\mathbb{Z})$$

$$(a_j) \mapsto \sum_{j=0}^{m-1} a_j$$

We define a natural action of $N$ on $Y$. For any $g = (a_0, \cdots, a_{m-1}) \in N$, the action of $g$ on $Y$ is induced by

$$g \cdot y_j := \zeta_r^{a_j} y_j, \quad \forall\ 0 \le j \le m-1,$$

where recall $\zeta_r$ is a fixed $r$-th primitive root of unity.

**Proposition 2.3.** The following relations between $Y$ and $X$ hold:

(1) The map $\pi : Y \to \mathbb{P}^n$, $[y_0 : \cdots : y_{m-1}] \mapsto [y_0^r : \cdots : y_n^r]$ defines a Galois covering of degree $r^{m-1}$.
(2) $X \simeq Y/N_1$.
(3) Under the isomorphism in (2), we can identify $H^n(X, \mathbb{Q})$ with $H^n(Y, \mathbb{Q})^{N_1}$, where $H^n(Y, \mathbb{Q})^{N_1}$ denotes the subspace of invariants under $N_1$.

*Proof.* See Lemma 2.4, Proposition 2.5 and Proposition 2.6 in [11] for the $r = 2$ case. The same proof goes through in the general case without difficulty. Here we only indicate that (3) follows from (2) by a standard application of Leray spectral sequence.

□

2.3. **A special locus.** There is an interesting locus in $\mathfrak{M}_{AR}$ where the cyclic cover $X$ of $\mathbb{P}^n$ is determined by a cyclic cover of $\mathbb{P}^1$ branched at some points.

Note that there exists a natural Galois covering with Galois group $S_n$, the permutation group of $n$ letters:

$$\gamma : (\mathbb{P}^1)^n \to Sym^n(\mathbb{P}^1) = \mathbb{P}^n.$$

Here the identification attaches to a divisor of degree $n$ the ray of its equation in $H^0(\mathbb{P}^1, \mathcal{O}(n))$.

**Lemma 2.4.** Let $(p_1, \cdots, p_m)$ be a collection of $m$ distinct points on $\mathbb{P}^1$, and put $H_i = \gamma(\{p_i\} \times \mathbb{P}^1 \times \cdots \times \mathbb{P}^1)$. Then $(H_1, \cdots, H_m)$ is an arrangement of hyperplanes in general position.

*Proof.* The divisors of degree $n$ in $\mathbb{P}^1$ containing a given point form a hyperplane and, as a divisor of degree $n$ cannot contain $n + 1$ distinct points, no $n + 1$ hyperplanes in the arrangement do meet. □

**Remark 2.2.** Written in homogeneous coordinates, if $p_i = [a : b] \in \mathbb{P}^1$, then the corresponding $H_i = [a^n : a^{n-1}b : a^{n-2}b^2 : \cdots : b^n] \in \check{\mathbb{P}}^n$.

Given $p_i$, $H_i$ $(1 \leq i \leq m)$ as in Lemma 2.4, let $C$ be the degree $r$ cyclic cover of $\mathbb{P}^1$ branched at $p_1, \cdots, p_m$, and let $X_C$ be the degree $r$ cyclic cover of $\mathbb{P}^n$ branched along $H_1, \cdots, H_m$. Then the cyclic covering structures induce natural actions of the cyclic group $\mathbb{Z}/r\mathbb{Z}$ on $C$ and $X_C$.

The group $(\mathbb{Z}/r\mathbb{Z})^n$ and the permutation group $S_n$ act naturally on the product $C^n$. These actions induce an action of the semi-direct product $(\mathbb{Z}/r\mathbb{Z})^n \rtimes S_n$ on $C^n$. Let $N'$ be the kernel of the summation homomorphism:

$$(\mathbb{Z}/r\mathbb{Z})^n \to \mathbb{Z}/r\mathbb{Z}$$

$$(a_i) \mapsto \sum_{i=1}^{n} a_i$$

Then we have:

**Lemma 2.5.** There exists a natural isomorphism: $X_C \simeq C^n/N_2$, where $N_2 := N' \rtimes S_n$.

*Proof.* Essentially the same proof of Lemma 2.5 in [20]. □

**Proposition 2.6.** For each $1 \le i \le r-1$, we have

$$H^n(X_C, \mathbb{C})_{(i)} \simeq \wedge^n H^1(C, \mathbb{C})_{(i)}.$$

*Proof.* It follows from Lemma 2.5 and the Künneth formula. For more details about the $i = 1$ case, we refer to Proposition 2.7 in [20]. Other cases can be treated in a completely analogous way.

□

Now we consider the bilinear forms $Q$ on $H^n(X_C, \mathbb{R})$ and $Q_C$ on $H^1(C, \mathbb{R})$ induced by cup products, together with the associated Hermitian forms:

$$H(\alpha, \beta) := \sqrt{-1}^n Q(\alpha, \bar{\beta}); \quad H_C(\alpha, \beta) := \sqrt{-1} Q_C(\alpha, \bar{\beta})$$

By Proposition 2.3, $X_C$ is a quotient of a smooth projective variety, so the natural mixed Hodge structure on $X_C$ is pure (cf. [4]). The Hodge decomposition

$$H^n(X_C, \mathbb{C}) = \oplus_{p+q=n} H^{p,q}(X_C)$$

is compatible with the $\mathbb{Z}/r\mathbb{Z}$-action. Hence for each $1 \le i \le r-1$, we have a decomposition of the $i$-eigenspace:

$$H^n(X_C, \mathbb{C})_{(i)} = \oplus_{p+q=n} H^{p,q}(X_C)_{(i)}.$$

Similar decomposition

$$H^1(C, \mathbb{C})_{(i)} = H^{1,0}(C)_{(i)} \oplus H^{0,1}(C)_{(i)}$$

holds on $C$, and Proposition 2.6 implies

$$H^{p,q}(X_C)_{(i)} \simeq \wedge^p H^{1,0}(C)_{(i)} \otimes \wedge^q H^{0,1}(C)_{(i)}.$$

Since $H_C$ is positive definite on $H^{1,0}(C)$ and negative definite on $H^{0,1}(C)$, the correspondence between $H$ and $H_C$ implies there exists $\epsilon(n) = \pm 1$, such that $\epsilon(n)H$ is positive definite on

$$V_{\text{even}} := \bigoplus_{p+q=n,\ q \text{ even}} H^{p,q}(X_C)_{(i)},$$

and negative definite on

$$V_{\text{odd}} := \bigoplus_{p+q=n,\ q \text{ odd}} H^{p,q}(X_C)_{(i)}.$$

As $H^n(X_C, \mathbb{C})_{(i)} = V_{\text{even}} \oplus V_{\text{odd}}$, we see for each $1 \le i \le r-1$, the restriction of $\epsilon(n)H$ on $H^n(X_C, \mathbb{C})_{(i)}$ is non-degenerate, with signature $(p_i, q_i)$, and by the following lemma, $p_i$, $q_i$ can be written down explicitly:

$$p_i = \sum_{j=0}^{\llcorner \frac{n}{2} \lrcorner} \binom{m-ki-1}{n-2j} \binom{ki-1}{2j},$$

$$q_i = \sum_{j=0}^{\llcorner \frac{n-1}{2} \lrcorner} \binom{m-ki-1}{n-2j-1} \binom{ki-1}{2j+1},$$

where $k = \frac{m}{r}$.

**Lemma 2.7.** For $1 \leq i \leq r-1$, we have $dim H^{1,0}(C)_{(i)} = m - \frac{mi}{r} - 1$, $dim H^{0,1}(C)_{(i)} = \frac{mi}{r} - 1$.

*Proof.* See [16], (2.7). □

## 3. A Picard-Lefschetz Type Formula

In this section, we study the monodromy transformation of a one-parameter degeneration of cyclic covers of $\mathbb{P}^n$ branched along hyperplane arrangements. It turns out that this monodromy transformation is given by a form analogous to the classical Picard-Lefschetz formula. It is obtained by applying Pham's generalised Picard-Lefshcetz formula to the corresponding degeneration of Kummer covers.

3.1. **The localization principle.** If the central fiber of a one-parameter degeneration of varieties admits only finitely many isolated singularities, then the monodromy transformation of this family is determined by the local monodromies around these singularities. This localization principle allows us to reduce the monodromy computation to a local question about the isolated singularities.

Here is the general situation. Suppose $Y$ is a topological space and $B$, $X_i$ $(1 \leq i \leq m)$ are closed subspaces of $Y$ such that $X_i$ $(1 \leq i \leq m)$ are pairwise disjoint and $Y = \cup_{i=1}^m X_i \cup B$. Suppose $\Phi : Y \to Y$ is a homeomorphism which restricts to the identity map on $B$ and which maps $X_i$ to $X_i$, for each $1 \leq i \leq m$. It can be imagined that the information of the restrictions $\Phi_i = \Phi|_{X_i} : X_i \to X_i$ determines the induced homomorphism $\Phi_* : H_n(Y, \mathbb{Z}) \to H_n(Y, \mathbb{Z})$ on the homology group. More precisely, let $A_i = B \cap X_i$ and consider the variation homomorphism:

$$Var_n(\Phi_i) : H_n(X_i, A_i, \mathbb{Z}) \to H_n(X_i, \mathbb{Z})$$

where for each $[c] \in H_n(X_i, A_i, \mathbb{Z})$ represented by a relative cycle $c$, we define the homology class $Var_n(\Phi_i)([c]) \in H_n(X_i, \mathbb{Z})$ to be the class represented by $\Phi_i(c) - c$, which is well-defined since $\Phi_i$ restricts to the identity map on $A_i$. Assume furthermore that the natural homomorphism induced by inclusions

$$\phi : \oplus_{i=1}^m H_n(X_i, A_i, \mathbb{Z}) \to H_n(Y, B, \mathbb{Z})$$

is an isomorphism (this holds by the Excision Theorem if $B$ contains a closed subset $Z$ such that $\cup_{i=1}^m A_i \cap Z = \emptyset$ and the pair $(\cup_{i=1}^m X_i, \cup_{i=1}^m A_i)$ is a deformation retraction of $(Y - Z, B - Z)$). Under this assumption, we can show

**Proposition 3.1.** The homomorphism $\Phi_* - id : H_n(Y,\ \mathbb{Z}) \to H_n(Y,\ \mathbb{Z})$ is equal to the composition

$$\begin{array}{ccc} H_n(Y,\ \mathbb{Z}) \xrightarrow{\phi_1} H_n(Y, B,\ \mathbb{Z}) \xrightarrow[\sim]{\phi^{-1}} & \oplus_{i=1}^m H_n(X_i, A_i,\ \mathbb{Z}) & \\ & \downarrow{\scriptstyle \oplus_i Var_n(\Phi_i)} & \\ & \oplus_{i=1}^m H_n(X_i,\ \mathbb{Z}) & \xrightarrow{\phi_2} H_n(Y,\ \mathbb{Z}) \end{array}$$

where $\phi_1$ is the natural homomorphism and $\phi_2$ is induced by the inclusions.

*Proof.* This is standard and it follows directly from the definition of $Var_n(\Phi_i)$. See the discussion in [13], Chapter 3, for example.

□

3.2. **Restatement of Pham's results.** In [17], F. Pham considered the singularity defined by the polynomial equation $\sum_{i=1}^{n+1} z_i^{\mu_i} = 0$ and derived a generalised Picard-Lefschetz formula about the monodromy transformation around these singularities. Here $\mu_i$ are positive integers. In this subsection, we restate a special case of Pham's result in a form convenient for us. To be precise, for $r \geq 2$, consider the function $f : \mathbb{C}^{n+1} \to \mathbb{C}$ defined by $f(z) = \sum_{i=1}^{n+1} z_i^r$. It is known for each sufficiently small $\epsilon > 0$, we can pick a positive number $\eta \ll \epsilon$, such that by putting $\bar{\mathcal{X}} := \{z \in \mathbb{C}^{n+1} : |z| \leq \epsilon, |f(z)| \leq \eta\}$ and $S := \{t \in \mathbb{C} : |t| < \eta\}$, the map

$$f : \bar{\mathcal{X}} \to S$$

restricts to a locally trivial $C^\infty$-fibration of compact manifolds with boundary over $S^* = S - \{0\}$, and $f|_{\partial\bar{\mathcal{X}}} : \partial\bar{\mathcal{X}} \to S$ is a trivial $C^\infty$-fibration (cf. [15]). Consequently, for $\rho \in (0, \eta)$, the geometric monodromy $\Phi : \bar{\mathcal{X}}_\rho \xrightarrow{\sim} \bar{\mathcal{X}}_\rho$ of the fiber over $\rho$ restricts to the identity map on $\partial\bar{\mathcal{X}}_\rho$.

For ease of notations, we fix $\rho \in (0, \eta)$ and let $X = \bar{\mathcal{X}}_\rho$, $A = \partial\bar{\mathcal{X}}_\rho$. Next we want to compute the variation homomorphism

$$Var_n(\Phi) : H_n(X, A,\ \mathbb{Q}(\zeta_r)) \to H_n(X,\ \mathbb{Q}(\zeta_r)).$$

First note that since $(X, A)$ is a compact orientable manifold with boundary, Lefschetz duality identifies $H_n(X, A,\ \mathbb{Z})$ with $H^n(X,\ \mathbb{Z})$, and the natural paring

$$H^n(X,\ \mathbb{Z}) \times H_n(X,\ \mathbb{Z}) \to \mathbb{Z}$$

induces the intersection paring

$$Q : H_n(X, A,\ \mathbb{Z}) \times H_n(X,\ \mathbb{Z}) \to \mathbb{Z}.$$

We define

$$H : H_n(X, A,\ \mathbb{Q}(\zeta_r)) \times H_n(X,\ \mathbb{Q}(\zeta_r)) \to \mathbb{C}$$

by the formula $H(\alpha, \beta) := \sqrt{-1}^n Q(\alpha, \bar{\beta})$.

Let $G=(\mathbb{Z}/r\mathbb{Z})^{n+1}$. We define a $G$-action on $\mathbb{C}^{n+1}$ by

$$G\times\mathbb{C}^{n+1}\to\mathbb{C}^{n+1}$$
$$((a_i)_{i=1}^{n+1},(z_i)_{i=1}^{n+1})\mapsto(\zeta_r^{a_i}z_i)_{i=1}^{n+1}.$$

This induces a $G$-action on the fibration $f:\bar{\mathcal{X}}\to S$. Obviously $Var_n(\Phi)$ is $G$-equivariant, and $H_n(X,A,\ \mathbb{Q}(\zeta_r))$ admits a decomposition

$$H_n(X,A,\ \mathbb{Q}(\zeta_r))=\bigoplus_{(\mu_i)\in(\mathbb{Z}/r\mathbb{Z})^{n+1}}H_n(X,A,\ \mathbb{Q}(\zeta_r))_{(\mu_1,\cdots,\mu_{n+1})}$$

where $H_n(X,A,\ \mathbb{Q}(\zeta_r))_{(\mu_1,\cdots,\mu_{n+1})}$ is the following character subspace of $H_n(X,A,\ \mathbb{Q}(\zeta_r))$:

$$\{\alpha\in H_n(X,A,\ \mathbb{Q}(\zeta_r))|g_*(\alpha)=\zeta_r^{\sum_{i=1}^{n+1}a_i\mu_i}\alpha,\ \forall g=(a_i)\in G\}.$$

Similar decompositions hold:

$$H_n(X,\ \mathbb{Q}(\zeta_r))=\bigoplus_{(\mu_i)\in(\mathbb{Z}/r\mathbb{Z})^{n+1}}H_n(X,\ \mathbb{Q}(\zeta_r))_{(\mu_1,\cdots,\mu_{n+1})},$$

$$Var_n(\Phi)=\bigoplus_{(\mu_i)\in(\mathbb{Z}/r\mathbb{Z})^{n+1}}Var_n(\Phi)_{(\mu_1,\cdots,\mu_{n+1})},$$

where $Var_n(\Phi)_{(\mu_1,\cdots,\mu_{n+1})}:H_n(X,A,\ \mathbb{Q}(\zeta_r))_{(\mu_1,\cdots,\mu_{n+1})}\to H_n(X,\ \mathbb{Q}(\zeta_r))_{(\mu_1,\cdots,\mu_{n+1})}$ is the $(\mu_1,\cdots,\mu_{n+1})$-character part of $Var_n(\Phi)$.

**Proposition 3.2** (Pham)**.** The following statements hold:

(i)

$$dim_{\mathbb{Q}(\zeta_r)}H_n(X,\ \mathbb{Q}(\zeta_r))_{(\mu_1,\cdots,\mu_{n+1})}=\begin{cases}1, & \text{if }\ \mu_i\neq0\in\mathbb{Z}/r\mathbb{Z},\ \ \forall\ 1\leq i\leq n+1;\\ 0, & \text{otherwise .}\end{cases}$$

(ii)

$$dim_{\mathbb{Q}(\zeta_r)}H_n(X,A,\ \mathbb{Q}(\zeta_r))_{(\mu_1,\cdots,\mu_{n+1})}=\begin{cases}1, & \text{if }\ \mu_i\neq0\in\mathbb{Z}/r\mathbb{Z},\ \ \forall\ 1\leq i\leq n+1;\\ 0, & \text{otherwise .}\end{cases}$$

(iii) For each $(\mu_1,\cdots,\mu_{n+1})\in(\mathbb{Z}/r\mathbb{Z})^{n+1}$ with $\mu_i\neq0\in\mathbb{Z}/r\mathbb{Z},\ \forall\ 1\leq i\leq n+1$, there is a generator $e_{(\mu_1,\cdots,\mu_{n+1})}\in H_n(X,\ \mathbb{Q}(\zeta_r))_{(\mu_1,\cdots,\mu_{n+1})}$, such that $\forall\ \alpha\in H_n(X,A,\ \mathbb{Q}(\zeta_r))_{(\mu_1,\cdots,\mu_{n+1})}$, the variation homomorphism is given by the generalised Picard-Lefschetz formula:

$$Var_n(\Phi)_{(\mu_1,\cdots,\mu_{n+1})}(\alpha)=c_{(\mu_1,\cdots,\mu_{n+1})}H(\alpha,e_{(\mu_1,\cdots,\mu_{n+1})})e_{(\mu_1,\cdots,\mu_{n+1})},\tag{1}$$

where $c_{(\mu_1,\cdots,\mu_{n+1})}=\frac{-(-1)^{n(n+1)/2}r^{n+1}}{\sqrt{-1}^n\prod_{j=1}^{n+1}(1-\zeta_r^{-\mu_j})}\in\mathbb{Q}(\zeta_r,\sqrt{-1})^*$. Moreover, we have

$$H(e_{(\mu_1,\cdots,\mu_{n+1})},e_{(\mu_1,\cdots,\mu_{n+1})})=(\zeta_r^{\sum_{i=1}^{n+1}\mu_i}-1)c^{-1}_{(\mu_1,\cdots,\mu_{n+1})}.$$

*Proof.* (i): Let $\omega_i = (\delta_{ij})_{j=1}^{n+1} \in G = (\mathbb{Z}/r\mathbb{Z})^{n+1}$ where

$$\delta_{ij} = \begin{cases} 1, & \text{if } j = i; \\ 0, & \text{if } j \neq i. \end{cases}$$

Then Theorem 1 in Pham [17] asserts that as a $\mathbb{Z}[G]$-module, $H_n(X,\ \mathbb{Z})$ is isomorphic to the image of the $\mathbb{Z}[G]$-homomorphism

$$\mathbb{Z}[G] \xrightarrow{(1-\omega_1)(1-\omega_2)\cdots(1-\omega_{n+1})} \mathbb{Z}[G].$$

From this we can deduce (i) without difficulty.

(ii) follows from Corollary 1 in Pham *op. cit.* in a similar way.

(iii): Following Pham, given $\alpha \in H_n(X, A,\ \mathbb{Z})$ and $\beta \in H_n(X,\ \mathbb{Z})$, we define

$$(\alpha|\beta) := \sum_{g\in G} Q(\alpha,\ g_*\beta)g^{-1} \in \mathbb{Z}[G].$$

By the formula (Var I ) in Pham *op. cit.*, the variation homomorphism

$$Var_n(\Phi) : H_n(X, A,\ \mathbb{Z}) \to H_n(X,\ \mathbb{Z})$$

is given by

$$Var_n(\Phi)(\epsilon) = -(-1)^{n(n+1)/2} e, \tag{2}$$

where $\epsilon$ (resp. $e$ ) is a suitable $\mathbb{Z}[G]$-generator of $H_n(X, A,\ \mathbb{Z})$( resp. $H_n(X,\ \mathbb{Z})$), and

$$(\epsilon|e) = (1-\omega_1)(1-\omega_2)\cdots(1-\omega_{n+1}) \in \mathbb{Z}[G]. \tag{3}$$

In $H_n(X, A,\ \mathbb{Q}(\zeta_r))$, we can decompose $\epsilon$ as

$$\epsilon = \sum_{(\mu_1,\cdots,\mu_{n+1})\in(\mathbb{Z}/r\mathbb{Z})^{n+1}} \epsilon_{(\mu_1,\cdots,\mu_{n+1})} \tag{4}$$

in such a way that $\epsilon_{(\mu_1,\cdots,\mu_{n+1})}$ linearly spans $H_n(X, A,\ \mathbb{Q}(\zeta_r))_{(\mu_1,\cdots,\mu_{n+1})}$. Similarly we have the decomposition

$$e = \sum_{(\mu_1,\cdots,\mu_{n+1})\in(\mathbb{Z}/r\mathbb{Z})^{n+1}} e_{(\mu_1,\cdots,\mu_{n+1})}. \tag{5}$$

Note that $\bar{\epsilon}_{(\mu_1,\cdots,\mu_{n+1})} = \epsilon_{(-\mu_1,\cdots,-\mu_{n+1})}$ and $\bar{e}_{(\mu_1,\cdots,\mu_{n+1})} = e_{(-\mu_1,\cdots,-\mu_{n+1})}$.

We can obtain from (4) that

$$\epsilon_{(\mu_1,\cdots,\mu_{n+1})} = \frac{1}{r^{n+1}} \sum_{(a_1,\cdots,a_{n+1})\in(\mathbb{Z}/r\mathbb{Z})^{n+1}} \zeta_r^{-\sum_{i=1}^{n+1} a_i\mu_i} \omega_{1*}^{a_1}\omega_{2*}^{a_2}\cdots\omega_{n+1*}^{a_{n+1}}\epsilon. \tag{6}$$

Similarly, from (5), we get

$$e_{(\mu_1,\cdots,\mu_{n+1})} = \frac{1}{r^{n+1}} \sum_{(a_1,\cdots,a_{n+1})\in(\mathbb{Z}/r\mathbb{Z})^{n+1}} \zeta_r^{-\sum_{i=1}^{n+1} a_i\mu_i} \omega_{1*}^{a_1}\omega_{2*}^{a_2}\cdots\omega_{n+1*}^{a_{n+1}}e. \tag{7}$$

For each $g \in G$, $\alpha \in H_n(X, A,\ \mathbb{Z})$ and $\beta \in H_n(X,\ \mathbb{Z})$, it follows from the $G$-invariance of the intersection form $Q$ that

$$(g_*^{-1}\alpha|\beta) = (\alpha|g_*\beta) = g(\alpha|\beta).$$

Then plugging the expressions (6) and (7) into (3), we get

$$(\epsilon_{(\mu_1,\cdots,\mu_{n+1})}|e_{(-\mu_1,\cdots,-\mu_{n+1})}) = \frac{1}{r^{n+1}} \sum_{(a_1,\cdots,a_{n+1})\in(\mathbb{Z}/r\mathbb{Z})^{n+1}} \zeta_r^{\sum_{i=1}^{n+1} a_i\mu_i} \prod_{j=1}^{n+1} \omega_j^{a_j}(1-\omega_j)$$

From the definition of $(\cdot|\cdot)$, we can obtain the intersection number

$$Q(\epsilon_{(\mu_1,\cdots,\mu_{n+1})}, e_{(-\mu_1,\cdots,-\mu_{n+1})}) = \frac{1}{r^{n+1}} \prod_{j=1}^{n+1}(1-\zeta_r^{-\mu_j}). \tag{8}$$

Moreover, we know from (2) that

$$Var_n(\Phi)(\epsilon_{(\mu_1,\cdots,\mu_{n+1})}) = -(-1)^{n(n+1)/2} e_{(\mu_1,\cdots,\mu_{n+1})}.$$

Combining this expression with (8), we obtain finally that if $\mu_i \neq 0 \in \mathbb{Z}/r\mathbb{Z}$, $\forall\ 1 \leq i \leq n+1$, then

$$Var_n(\Phi)(\epsilon_{(\mu_1,\cdots,\mu_{n+1})}) = c_{(\mu_1,\cdots,\mu_{n+1})} H(\epsilon_{(\mu_1,\cdots,\mu_{n+1})}, e_{(\mu_1,\cdots,\mu_{n+1})}) e_{(\mu_1,\cdots,\mu_{n+1})}.$$

Since $\epsilon_{(\mu_1,\cdots,\mu_{n+1})}$ spans the $\mathbb{Q}(\zeta_r)$-linear space $H_n(X, A,\ \mathbb{Q}(\zeta_r))_{(\mu_1,\cdots,\mu_{n+1})}$, we get the generalised Picard-Lefschetz formula (1). Combining (7) and the formula

$$(e|e) = (-1)^{n(n+1)/2}(1-\omega_1)\cdots(1-\omega_{n+1})(1-\omega_1^{-1}\omega_2^{-1}\cdots\omega_{n+1}^{-1})$$

in Pham *op. cit.*, p. 340, we get the formula

$$H(e_{(\mu_1,\cdots,\mu_{n+1})}, e_{(\mu_1,\cdots,\mu_{n+1})}) = (\zeta_r^{\sum_{i=1}^{n+1}\mu_i} - 1)c_{(\mu_1,\cdots,\mu_{n+1})}^{-1}.$$

□

3.3. **The Picard-Lefschetz type formula for cyclic covers of $\mathbb{P}^n$.** Given complex numbers $a_{ij}$ $(1 \leq i \leq n,\ 1 \leq j \leq m-n-2)$, for each $t \in \mathbb{C}$, consider the hyperplane arrangement $\mathfrak{A}_t = (H_{1t}, \cdots, H_{mt})$ of $\mathbb{P}^n$ represented by the following $(n+1)\times m$ matrix:

$$B = (b_{ij}) = \begin{pmatrix} 1 & 0 & \cdots & 0 & t & 1 & \cdots & 1 \\ 0 & 1 & & 0 & -1 & a_{11} & \cdots & a_{1,m-n-2} \\ \vdots & & \ddots & \vdots & \vdots & \vdots & & \vdots \\ 0 & & & 1 & -1 & a_{n1} & \cdots & a_{n,m-n-2} \end{pmatrix}$$

Here the $j$-th column corresponds to the defining equation

$$\sum_{i=0}^{n} b_{ij}x_i = 0$$

of the hyperplane $H_{jt}$, and $[x_0 : \cdots : x_n]$ are the homogeneous coordinates on $\mathbb{P}^n$. We assume the hyperplane arrangement $(H_{1t}, \cdots, H_{n+1,t}, H_{n+3,t}, \cdots, H_{mt})$ is in general position. Moreover, we assume there exists a positive number $\epsilon > 0$, such that for each $t$ in the punctured disk $\Delta_\epsilon^* := \{t \in \mathbb{C} : 0 < |t| < \epsilon\}$, the hyperplane arrangement $\mathfrak{A}_t$ is in general position, and for $t = 0$, the divisor $\sum_{i=1}^m H_{i0}$ has simple normal crossings on $\mathbb{P}^n - \{[1 : 0 : \cdots : 0]\}$.

The one parameter family of hyperplane arrangements $\mathfrak{A}_t$ determines a family $\mathcal{X}_\epsilon^* \to \Delta_\epsilon^*$, in such a way that the fiber $\mathcal{X}_{\epsilon,t}^*$ over $t \in \Delta_\epsilon^*$ is the degree $r$ cyclic cover of $\mathbb{P}^n$ branched along $\mathfrak{A}_t$. Fixing a base point $\rho \in \Delta_\epsilon^*$, we are interested in the monodromy transformation on $H^n(\mathcal{X}_{\epsilon,\rho}^*, \mathbb{C})$.

Now consider the family of the Kummer covers $\mathcal{Y}_\epsilon \xrightarrow{\pi} \Delta_\epsilon$, where $\mathcal{Y}_\epsilon$ is the subvariety of $\mathbb{P}^{m-1} \times \Delta_\epsilon$ defined by the equations

$$ty_0^r - \sum_{j=1}^{n+1} y_j^r = 0;$$

$$y_{n+i+1}^r - (y_0^r + \sum_{j=1}^{n} a_{ji} y_j^r) = 0,\ 1 \le i \le m-n-2.$$

and $\pi$ is induced from the natural projection morphism $\mathbb{P}^{m-1} \times \Delta_\epsilon \to \Delta_\epsilon$. Here $[y_0 : \cdots : y_{m-1}]$ are the homogeneous coordinates on $\mathbb{P}^{m-1}$ and $\Delta_\epsilon := \{t \in \mathbb{C} : |t| < \epsilon\}$. Recall from subsection 2.2 that the group $N = \oplus_{j=0}^{m-1} \mathbb{Z}/r\mathbb{Z}$ acts on $\mathcal{Y}_\epsilon$ fiberwisely and we have the isomorphism $\mathcal{X}_\epsilon^* \simeq \pi^{-1}(\Delta_\epsilon^*)/N_1$, with $N_1 = Ker(N \xrightarrow{\sum} \mathbb{Z}/r\mathbb{Z})$. We first analyse the monodromy of the family $\pi$.

Note the critical locus $C_\pi$ of $\pi$ (i.e. the set of points of $\mathcal{Y}_\epsilon$ which are singular or where $\pi$ is not a submersion) is

$$\{([y_0 : \cdots : y_{m-1}], t) | y_0 = 1, t = y_1 = \cdots = y_{n+1} = 0, y_{n+2}^r = \cdots = y_{m-1}^r = 1\}.$$

For a set of elements $a_j \in \mathbb{Z}/r\mathbb{Z}$, $n+2 \le j \le m-1$, we use $p^{(a_{n+2},\cdots,a_{m-1})}$ to denote the point in $C_\pi$ with $y_j = \zeta_r^{a_j}$ $(n+2 \le j \le m-1)$. It can be seen by a direct computation that, each $p^{(a_{n+2},\cdots,a_{m-1})}$ admits an open neighborhood $U^{(a_{n+2},\cdots,a_{m-1})}$ in $\mathcal{Y}_\epsilon$ with the following commutative diagram (shrinking $\epsilon$ if necessary):

$$\begin{array}{ccc} U^{(a_{n+2},\cdots,a_{m-1})} & \xrightarrow[\sim]{\phi} & B_\eta \\ & {\scriptstyle\pi}\searrow \quad \swarrow{\scriptstyle f} & \\ & \Delta_\epsilon & \end{array} \tag{9}$$

where $B_\eta := \{z \in \mathbb{C}^{n+1} : |z| < \eta, |f(z)| < \epsilon\}$, $f(z) = z_1^r + \cdots + z_{n+1}^r$, and $\phi$ is given by

$$\phi([y_0 : \cdots : y_{m-1}], t) = (\frac{y_1}{y_0}, \cdots, \frac{y_{n+1}}{y_0}).$$

For the base point $\rho \in \Delta_\epsilon^*$, we put $Y = \pi^{-1}(\rho)$. The geometric monodromy $\tilde{\Phi} : Y \to Y$ of the family $\pi$ induces the monodromy transformation on the middle homology

$$\tilde{\Phi}_* : H_n(Y,\ \mathbb{Q}(\zeta_r)) \to H_n(Y,\ \mathbb{Q}(\zeta_r)).$$

As usual, the intersection form $Q$ induces the Hermitian form $H(\alpha,\beta) := \sqrt{-1}^n Q(\alpha,\bar{\beta})$ on $H_n(Y,\ \mathbb{Q}(\zeta_r))$. Define the set

$$M := \{(a_j) \in \oplus_{j=1}^{m-1}\mathbb{Z}/r\mathbb{Z} : a_j \neq 0,\ \forall\ 1 \leq j \leq n+1\}.$$

Then $\tilde{\Phi}_*$ can be represented as follows:

**Proposition 3.3.** For each $(a_j) \in M$, there exists $e_{(a_1,\cdots,a_{n+1})}^{(a_{n+2},\cdots,a_{m-1})} \in H_n(Y,\ \mathbb{Q}(\zeta_r))$, such that $\forall\ \alpha \in H_n(Y,\ \mathbb{Q}(\zeta_r))$, we have

$$\tilde{\Phi}_*(\alpha) - \alpha = \sum_{(a_j)\in M} c_{(a_1,\cdots,a_{n+1})} H(\alpha, e_{(a_1,\cdots,a_{n+1})}^{(a_{n+2},\cdots,a_{m-1})}) e_{(a_1,\cdots,a_{n+1})}^{(a_{n+2},\cdots,a_{m-1})},$$

where $c_{(a_1,\cdots,a_{n+1})} = \frac{-(-1)^{n(n+1)/2} r^{n+1}}{\sqrt{-1}^n \prod_{j=1}^{n+1}(1-\zeta_r^{-a_j})}$ is nonzero. Moreover, each $g = (b_j) \in N = \oplus_{j=0}^{m-1}\mathbb{Z}/r\mathbb{Z}$ acts on $e_{(a_1,\cdots,a_{n+1})}^{(a_{n+2},\cdots,a_{m-1})}$ in the following way:

$$g_* e_{(a_1,\cdots,a_{n+1})}^{(a_{n+2},\cdots,a_{m-1})} = \zeta_r^{\sum_{j=1}^{n+1} a_j(b_j-b_0)} e_{(a_1,\cdots,a_{n+1})}^{(a_{n+2}+b_{n+2}-b_0,\cdots,a_{m-1}+b_{m-1}-b_0)}.$$

*Proof.* By the definition of the action of $N$ on the family $\mathcal{Y}_\epsilon \xrightarrow{\pi} \Delta_\epsilon$, for each $g = (b_j) \in N$ and $p^{(a_{n+2},\cdots,a_{m-1})} \in C_\pi$, we have

$$g \cdot p^{(a_{n+2},\cdots,a_{m-1})} = p^{(a_{n+2}+b_{n+2}-b_0,\cdots,a_{m-1}+b_{m-1}-b_0)}.$$

As a consequence, for each $p^{(a_{n+2},\cdots,a_{m-1})} \in C_\pi$, we can shrink its open neighborhood $U^{(a_{n+2},\cdots,a_{m-1})}$ in the commutative diagram (9) so that it satisfies:

$$g \cdot U^{(a_{n+2},\cdots,a_{m-1})} = U^{(a_{n+2}+b_{n+2}-b_0,\cdots,a_{m-1}+b_{m-1}-b_0)}.$$

Then a combination of Proposition 3.1 and Proposition 3.2 completes the proof.

□

Let $X = \mathcal{X}_{\epsilon,\rho}^*$ be the fiber of $\mathcal{X}_\epsilon^*$ over the base point $\rho \in \Delta_\epsilon^*$. Then $X = Y/N_1$ and we have the Hermitian form $H$ on $H_n(X,\ \mathbb{Q}(\zeta_r))$ induced by the intersection pairing. The $\mathbb{Z}/r\mathbb{Z}$-action on $H_n(X,\ \mathbb{Q}(\zeta_r))$ induces a decomposition

$$H_n(X,\ \mathbb{Q}(\zeta_r)) = \oplus_{i=0}^{r-1} H_n(X,\ \mathbb{Q}(\zeta_r))_{(i)}.$$

Our main result in this section is the following

**Proposition 3.4.** The geometric monodromy $\Phi : X \to X$ induces the monodromy transformation $\Phi_* : H_n(X,\ \mathbb{Q}(\zeta_r)) \to H_n(X,\ \mathbb{Q}(\zeta_r))$, and $\forall\ 1 \leq i \leq r-1$, there exists a cycle $e_{(i)} \in H_n(X,\ \mathbb{Q}(\zeta_r))_{(i)}$, such that: $\forall\ \alpha \in H_n(X,\ \mathbb{Q}(\zeta_r))_{(i)}$, we have the formula

$$\Phi_*(\alpha) = \alpha + cH(\alpha, e_{(i)})e_{(i)},$$

where $c = \frac{-(-1)^{n(n+1)/2} r^{n+m}}{\sqrt{-1}^n (1-\zeta_r^{-i})^{n+1}} \in \mathbb{Q}(\zeta_r, \sqrt{-1})^*$ is a nonzero constant. Moreover, we have

$$H(e_{(i)}, e_{(i)}) = (\zeta_r^{(n+1)i} - 1)c^{-1}.$$

*Proof.* Since $H^n(X,\ \mathbb{Q}(\zeta_r)) = H^n(Y,\ \mathbb{Q}(\zeta_r))^{N_1}$ by Proposition (2.3), we can identify $H_n(X,\ \mathbb{Q}(\zeta_r))$ with the quotient space $H_n(Y,\ \mathbb{Q}(\zeta_r))/U$, where $U$ is the $\mathbb{Q}(\zeta_r)$-linear subspace of $H_n(Y,\ \mathbb{Q}(\zeta_r))$ spanned by elements in the set

$$\{g_*\alpha - \alpha : g \in N_1, \alpha \in H_n(Y,\ \mathbb{Q}(\zeta_r))\}.$$

With this identification, we can see that all the elements $e^{(a_{n+2},\cdots,a_{m-1})}_{(a_1,\cdots,a_{n+1})}$ in Proposition 3.3 become zero in $H_n(X,\ \mathbb{Q}(\zeta_r))$ except those with $a_1 = \cdots = a_{n+1} \in \mathbb{Z}/r\mathbb{Z}$. Moreover, each $e^{(a_{n+2},\cdots,a_{m-1})}_{(a_1,\cdots,a_{n+1})}$ can be identified with $e^{(0,\cdots,0)}_{(a_1,\cdots,a_{n+1})}$ in $H_n(X,\ \mathbb{Q}(\zeta_r))$ through the $N_1$-action. For $1 \le i \le r-1$, let $e_{(i)} = e^{(0,\cdots,0)}_{(i,\cdots,i)}$. Then the formula we want to prove follows from Proposition 3.3, by noticing the intersection form on $X$ differs a $r^{m-1}$ multiplication with that on $Y$.

□

For later use, we translate the preceding proposition (by Poincaré duality) to a formula about the cohomology groups.

**Proposition 3.5.** For each $1 \le i \le r-1$, there exists a cocycle $e_{(i)} \in H^n(X,\ \mathbb{Q}(\zeta_r))_{(i)}$, such that the monodromy transformation

$$\Phi_* \colon H^n(X, \mathbb{Q}(\zeta_r))_{(i)} \longrightarrow H^n(X, \mathbb{Q}(\zeta_r))_{(i)}$$

is given by the formula

$$\Phi_*(\alpha) = \alpha + cH(\alpha, e_{(i)})e_{(i)}, \qquad \forall\, \alpha \in H^n(X, \mathbb{Q}(\zeta_r))_{(i)},$$

where $c = \frac{-(-1)^{n(n+1)/2} r^{n+m}}{\sqrt{-1}^n (1-\zeta_r^{-i})^{n+1}} \in \mathbb{Q}(\zeta_r, \sqrt{-1})^*$ is a nonzero constant, and $H$ is the Hermitian form induced from the cup product paring, as defined in subsection 2.1. Moreover, we have

$$H(e_{(i)}, e_{(i)}) = (\zeta_r^{(n+1)i} - 1)c^{-1}.$$

## 4. Proof of the Main Result

In this section, we give the proof of our main result Theorem 2.1. We first introduce some moduli spaces related to $\mathfrak{M}_{AR}$ for technical reasons, and then as the algebraic preliminaries, we state our generalisation of Deligne's criterion and an alternative type result for semi-simple algebraic groups. After doing some reductions, we first give the proof of the double cover case Theorem 1.1, and then the proof of Theorem 2.1.

4.1. **Various moduli spaces.** Recall $\mathfrak{M}_{AR} \subset (\mathbb{P}^n)^m/PGL(n+1)$ is the coarse moduli space of ordered $m$ hyperplane arrangements in $\mathbb{P}^n$ in general position. We have the natural projections

$$(\mathbb{C}^{n+1} - \{0\})^m \to (\mathbb{P}^n)^m \to (\mathbb{P}^n)^m/PGL(n+1).$$

Let $\mathfrak{M}' \subset (\mathbb{P}^n)^m$ and $\widetilde{\mathfrak{M}} \subset (\mathbb{C}^{n+1}-\{0\})^m$ be the inverse images of $\mathfrak{M}_{AR}$ under the above projections. Then it is easy to see the projection $\widetilde{\mathfrak{M}} \to \mathfrak{M}'$ is a $(\mathbb{C}^*)^m$-bundle, and $\mathfrak{M}' \to \mathfrak{M}_{AR}$ is a trivial $PGL(n+1)$-bundle.

Recall we have constructed the universal family $\mathcal{X}_{AR} \xrightarrow{f} \mathfrak{M}_{AR}$ of degree $r$ cyclic covers of $\mathbb{P}^n$ branched along hyperplane arrangements in general position. In a similar way, we can construct a universal family $\widetilde{\mathcal{X}} \xrightarrow{\tilde{f}} \widetilde{\mathfrak{M}}$ such that $\forall\ \tilde{s} = (L_1, \cdots, L_m) \in \widetilde{\mathfrak{M}}$, the fiber $\tilde{f}^{-1}(\tilde{s})$ is the degree $r$ cyclic cover of $\mathbb{P}^n$ branched along the hyperplane arrangement given by the zero sets of the linear forms $L_1, \cdots, L_m$. We fix a base point $\tilde{s} \in \widetilde{\mathfrak{M}}$ lying above the previously chosen base point $s$ of $\mathfrak{M}_{AR}$, and identify the fiber $\tilde{f}^{-1}(\tilde{s})$ with $X = f^{-1}(s)$.

In the $n=1$ case, to emphasize the relation with curves, we use $\widetilde{\mathfrak{M}}_C$ (resp. $\widetilde{\mathcal{C}}$) to denote $\widetilde{\mathfrak{M}}$ (resp. $\widetilde{\mathcal{X}}$). So each fiber of the family $\widetilde{\mathcal{C}} \to \widetilde{\mathfrak{M}}_C$ is a degree $r$ cyclic cover of $\mathbb{P}^1$ branched along $m$ distinct points.

By the constructions in subsection 2.3, we have an embedding $\widetilde{\mathfrak{M}}_C \subset \widetilde{\mathfrak{M}}$, and the whole picture can be summarized as the following diagram:

$$\begin{array}{ccccccc} \widetilde{\mathcal{C}} & & \widetilde{\mathcal{X}} & & & & \mathcal{X}_{AR} \\ \downarrow & & \downarrow{\scriptstyle \tilde{f}} & & & & \downarrow{\scriptstyle f} \\ \widetilde{\mathfrak{M}}_C & \hookrightarrow & \widetilde{\mathfrak{M}} & \longrightarrow & \mathfrak{M}' & \longrightarrow & \mathfrak{M}_{AR} \end{array} \tag{10}$$

Suppose the base point $\tilde{s} \in \widetilde{\mathfrak{M}}$ belongs to $\widetilde{\mathfrak{M}}_C$, and let $C = \widetilde{\mathcal{C}}_{\tilde{s}}$ be the fiber of $\widetilde{\mathcal{C}}$ over $\tilde{s}$. For each integer $1 \leq i \leq \llcorner \frac{r}{2} \lrcorner$, We define $\widetilde{Mon}_{C,(i)}$ to be the smallest real algebraic subgroup of $Aut(H^1(C,\ \mathbb{C})_{(i)})$ containing the image of the monodromy representation

$$\tilde{\rho}_{C,i} : \pi_1(\widetilde{\mathfrak{M}}_C, \tilde{s}) \to Aut(H^1(C,\ \mathbb{C})_{(i)})$$

and define $\widetilde{Mon}_{(i)}$ to be the smallest real algebraic subgroup of $Aut(H^n(X,\ \mathbb{C})_{(i)})$ containing the image of the monodromy representation

$$\tilde{\rho}_i : \pi_1(\widetilde{\mathfrak{M}}, \tilde{s}) \to Aut(H^n(X,\ \mathbb{C})_{(i)}).$$

By Proposition 2.6, the isomorphism $H^n(X,\ \mathbb{C})_{(i)} \simeq \wedge^n H^1(C,\ \mathbb{C})_{(i)}$ implies the following commutative diagram:

$$\begin{array}{ccc} Aut(\wedge^n H^1(C,\ \mathbb{C})_{(i)}) & \xrightarrow{\sim} & Aut(H^n(X,\ \mathbb{C})_{(i)}) \\ \uparrow & & \uparrow \\ \cup & & \cup \\ \widetilde{Mon}_{C,(i)} & \longrightarrow & \widetilde{Mon}_{(i)} \end{array} \tag{11}$$

We can determine $\widetilde{Mon}^0_{C,(i)}$ by the following proposition.

**Proposition 4.1.** (1) If $r=2$, then

$$\widetilde{Mon}^0_{C,(1)} = Aut(H^1(C,\ \mathbb{R}),Q) = Sp(H^1(C,\ \mathbb{R}),Q);$$

(2) If $r>2$, $1\le i<\frac{r}{2}$ and $mi\ge 2r$, then

$$\widetilde{Mon}^0_{C,(i)} = SU(H^1(C,\ \mathbb{C})_{(i)},H).$$

*Proof.* (1) follows from [1] and (2) follows from Theorem 5.1.1 in [18]. □

Next we study the vanishing cycles on $X$. Note the discriminant locus $\Delta := (\mathbb{C}^{n+1}-\{0\})^m\backslash\widetilde{\mathfrak{M}}$ is a reducible hypersurface and its irreducible decomposition can be written as

$$\Delta = \sum_{1\le k_1<k_2<\cdots<k_{n+1}\le m} \Delta_{k_1k_2\cdots k_{n+1}}$$

where $\Delta_{k_1k_2\cdots k_{n+1}}$ is the collections of $(L_1,\cdots,L_m)\in(\mathbb{C}^{n+1}-\{0\})^m$ such that the hyperplanes defined by the zeros of the linear forms $L_{k_i}(1\le i\le n+1)$ have nonempty intersections in $\mathbb{P}^n$.

For each smooth point $p$ on $\Delta_{k_1k_2\cdots k_{n+1}}$, we can take a small open disk $D$ centered at $p$ and intersecting $\Delta_{k_1k_2\cdots k_{n+1}}$ transversally at $p$. Let $z$ be the local coordinate on $D$ with center $p$. For a small positive number $\epsilon$, we define a loop $u$ in $\widetilde{\mathfrak{M}}$ by

$$t\mapsto z=\epsilon\cdot e^{2\pi\sqrt{-1}t}$$

Taking a path $v$ in $\widetilde{\mathfrak{M}}$ from the base point $\tilde{s}$ to $u(0)=u(1)$, we define a loop

$$\gamma_{k_1k_2\cdots k_{n+1}} = vuv^{-1}.$$

This is a loop with the base point $\tilde{s}$, and we call any loop of such type a meridian of $\Delta_{k_1k_2\cdots k_{n+1}}$, following Carlson and Toledo [3]. By the discussion in subsection 3.3, there exists a vanishing cycle $e_{k_1k_2\cdots k_{n+1},(i)}\in H^n(X,\ \mathbb{C})_{(i)}$, such that the monodromy $\tilde{\rho}_i(\gamma_{k_1k_2\cdots k_{n+1}}): H^n(X,\ \mathbb{C})_{(i)}\to H^n(X,\ \mathbb{C})_{(i)}$ is given by the Picard-Lefschetz type formula

$$\tilde{\rho}_i(\gamma_{k_1k_2\cdots k_{n+1}})(\alpha) = \alpha + cH(\alpha, e_{k_1k_2\cdots k_{n+1},(i)})e_{k_1k_2\cdots k_{n+1},(i)}. \tag{12}$$

Here $c=\frac{-(-1)^{n(n+1)/2}r^{n+m}}{\sqrt{-1}^n(1-\zeta_r^{-i})^{n+1}}\in\mathbb{Q}(\zeta_r,\sqrt{-1})^*$ is a nonzero constant, and

$$H(e_{k_1k_2\cdots k_{n+1},(i)}, e_{k_1k_2\cdots k_{n+1},(i)}) = (\zeta_r^{(n+1)i}-1)c^{-1}.$$

Define the set

$$R := \{\gamma\cdot e_{k_1k_2\cdots k_{n+1},(i)} | \gamma\in\widetilde{Mon}_{(i)},\ 1\le k_1<k_2<\cdots<k_{n+1}\le m\}$$

which is a finite union of $\widetilde{Mon}_{(i)}$-orbits in $H^n(X,\ \mathbb{C})_{(i)}$.

**Proposition 4.2.** Suppose $1\le i\le r-1$ and $mi\ge 2r$, then the set $R$ spans the $\mathbb{C}$-linear space $H^n(X,\ \mathbb{C})_{(i)}$.

*Proof.* [1]Consider the set

$$S := \{\gamma \circ \gamma_{k_1 k_2 \cdots k_{n+1}} \circ \gamma^{-1} | \gamma \in \pi_1(\widetilde{\mathfrak{M}}, \tilde{s}), 1 \leq k_1 < k_2 < \cdots < k_{n+1} \leq m\}.$$

Since $(\mathbb{C}^{n+1} \backslash \{0\})^m$ is simply connected, it is well known that $S$ generates $\pi_1(\widetilde{\mathfrak{M}}, \tilde{s})$.

Suppose $R$ does not span $H^n(X, \mathbb{C})_{(i)}$. Since $H$ is non-degenerate on $H^n(X, \mathbb{C})_{(i)}$, we can find a nonzero element $\alpha_0 \in H^n(X, \mathbb{C})_{(i)}$ such that $H(\alpha_0, e) = 0$, for each $e \in R$. Then the Picard-Lefschetz type formula (12) implies that $\alpha_0$ is $\pi_1(\widetilde{\mathfrak{M}}, \tilde{s})$-invariant, and hence is $\widetilde{Mon}_{(i)}$-invariant. On the other hand, by the commutative diagram (11), we can view $\alpha_0$ as an element in $\wedge^n H^1(C, \mathbb{C})_{(i)}$ and it is $\widetilde{Mon}^0_{C,(i)}$-invariant. But this contradicts Proposition 4.1 since for $Sp(H^1(C, \mathbb{R})_{(1)})$ (resp. $SU(H^1(C, \mathbb{C})_{(i)})$), the representation on $\wedge^n H^1(C, \mathbb{R})_{(1)}$ (resp. $\wedge^n H^1(C, \mathbb{C})_{(i)}$) does not admit any nontrivial one dimensional invariant subspace. □

4.2. **A generalisation of Deligne's criterion.** In this subsection, we modify Deligne's criterion of Zariski density in three cases. The first two cases are from Deligne's original versions, and the last case is from Carlson and Toledo's complex reflection version in [3]. Rather than repeat the proofs from [6] or [3] word by word with obvious modifications, we remark that in the first two cases, our condition"$R$ consists of a finite union of $M$-orbits" guarantees that $R$ is a constructible subset of $V$, and hence it contains a Zariski dense open subset in its closure. In the last case, the same condition guarantees that $R$ is a semialgebraic subset of the underlying real linear space $V$, and this is exactly what Carlson and Toledo need in their proof.

**Proposition 4.3** (cf. Deligne [6], section 4.4)**.** Given a finite dimensional $\mathbb{C}$-linear space $V$ equipped with a non-degenerate alternating form $Q$, and a subset $R \subset V$ spanning $V$. Let $M$ be the smallest algebraic subgroup $M$ of $Sp(V, Q)$ containing the transvections $x \mapsto x + Q(x, \delta)\delta \ \ (\delta \in R)$. Suppose $R$ consists of a finite union of $M$-orbits. Then one of the followings holds:

- There exists a nontrivial $M$-invariant subspace $U$ of $V$ such that $R \subset U \cup U^{\perp}$.
- $M = Sp(V, Q)$.

**Proposition 4.4** (cf. Deligne [6], section 4.4)**.** Given a finite dimensional $\mathbb{C}$-linear space $V$ equipped with a non-degenerate symmetric form $Q$, and a subset $R \subset V$ spanning $V$. Assume the elements of $R$ satisfy $Q(\delta, \delta) = 2$ and let $M$ be the smallest algebraic subgroup $M$ of $O(V, Q)$ containing the reflections $x \mapsto x - Q(x, \delta)\delta \ \ (\delta \in R)$. Suppose $R$ consists of a finite union of $M$-orbits. Then one of the followings holds:

[1]**Correction notice.** The proof of Proposition 4.2 reproduced in this author manuscript contains an incorrect representation-theoretic assertion. A corrected proof, with a strengthened spanning statement, is given in the appended corrigendum.

- There exists a nontrivial $M$-invariant subspace $U$ of $V$ such that $R \subset U \cup U^{\perp}$.
- $M$ is finite.
- $M = O(V, Q)$.

**Proposition 4.5** (cf. Carlson and Toledo [3], Theorem 7.2)**.** Let $\epsilon = \pm 1$ be fixed. Given a finite dimensional $\mathbb{C}$-linear space $V$ equipped with a non-degenerate Hermitian form $h$, and a subset $R \subset V$ spanning $V$. Fix a root of unity $\lambda \neq \pm 1$. Assume the elements of $R$ satisfy $h(\delta, \delta) = \epsilon$ and let $M$ be the smallest algebraic subgroup $M$ of $U(V, h)$ containing the complex reflections $x \mapsto x + \epsilon(\lambda - 1)h(x, \delta)\delta \quad (\delta \in R)$. Suppose the signature $(p, q)$ of $h$ satisfies $p + q > 1$, and that $R$ consists of a finite union of $M$-orbits. Then one of the followings holds:

- There exists a nontrivial $M$-invariant subspace $U$ of $V$ such that $R \subset U \cup U^{\perp}$.
- $M$ is finite.
- The adjoint group $PM = PU(V, h)$.

4.3. **An alternative type result.** Given a finite dimensional complex vector space $W$ and suppose $1 \leq p \leq dimW - 1$, then the complex algebraic group $SL(W)$ acts naturally on the $p$-th wedge product $V := \wedge^p W$. In this way, we get an embedding $SL(W) \hookrightarrow SL(V)$. We call an element $g \neq Id \in SL(V)$ is a *pseudo-reflection* on $V$ if $rank(g - Id) = 1$, i.e., if there is a $g$-invariant and codimension one subspace $V_1$ of $V$ such that $g$ restricts to the identity map on $V_1$. We have the following alternative type result:

**Proposition 4.6.** Under the above embedding $SL(W) \hookrightarrow SL(V)$, if $G$ is a complex semi-simple connected algebraic group lying between $SL(W)$ and $SL(V)$, then $G$ is either $SL(W)$ or $SL(V)$. If moreover $G$ contains a pseudo-reflection on $V$, then $G = SL(V)$.

*Proof.* If $dimW = p+1$, then $V \simeq W$ and the embedding $SL(W) \hookrightarrow SL(V)$ is an isomorphism. Hence in this case our assertions hold trivially.

If $dimW \geq p + 2$, the first assertion follows from the Lie algebra version Proposition 6.9 in [20] (note the representation of $G$ on $V$ is automatically irreducible, as $G$ contains $SL(W)$ and $SL(W)$ acts irreducibly on $V = \wedge^p W$).

To prove the second assertion, let $\varphi_1 : SL(W) \hookrightarrow G$ and $\varphi_2 : G \hookrightarrow SL(V)$ be the embeddings. Suppose $\varphi_1$ is an isomorphism and $g \in G$ is a pseudo-reflection on $V$, then there exists $h \in SL(W)$ such that $\varphi_1(h) = g$. Since $g$ is obviously an unipotent element in $SL(V)$, we deduce $h$ is an unipotent element in $SL(W)$. Then explicit matrix computations show that $rank(\varphi_2(\varphi_1(h)) - Id) \geq 2$. This contradicts the assumption that $g = \varphi_1(h)$ is a pseudo-reflection on $V$. So we get $G = SL(V)$. □

4.4. **Two reductions.** Before starting the proof of the main result, we do some reductions.

**Claim 4.1.** Theorem 1.1 implies Theorem 2.1, (1).

Given a hyperplane arrangement $(H_1, \cdots, H_m)$ in $\mathbb{P}^n$ in general position, recall we have the degree $r$ cyclic cover $X \xrightarrow{\pi} \mathbb{P}^n$ branched along $\sum_{j=1}^m H_j$. Suppose $r$ admits a nontrivial factor $i$ and writing $r^{'} = \frac{r}{i}$, we can construct the degree $r^{'}$ cyclic cover $X^{'} \xrightarrow{\pi^{'}} \mathbb{P}^n$ branched along $\sum_{j=1}^m H_j$, and clearly there exists a natural morphism $\phi$ making the following diagram commutative

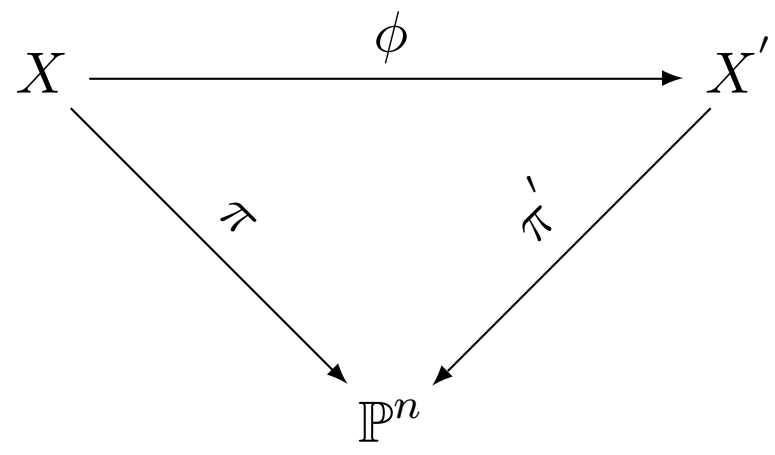

We have the $\mathbb{Z}/r\mathbb{Z}$-eigen space decomposition $H^n(X,\ \mathbb{C}) = \oplus_{j=0}^{r-1} H^n(X,\ \mathbb{C})_{(j)}$ and the $\mathbb{Z}/r^{'}\mathbb{Z}$-eigen space decomposition $H^n(X^{'},\ \mathbb{C}) = \oplus_{j=0}^{r^{'}-1} H^n(X^{'},\ \mathbb{C})_{(j)}$. Then it is not difficult to see that by taking $i = \frac{r}{2}$, our claim follows from the next proposition.

**Proposition 4.7.** $\phi$ induces an isomorphism $\phi_* : H^n(X^{'},\ \mathbb{C})_{(1)} \xrightarrow{\sim} H^n(X,\ \mathbb{C})_{(i)}$.

*Proof.* Let $\mathfrak{F} := \pi_*\mathbb{C}$ and $\mathfrak{F}^{'} := \pi^{'}_*\mathbb{C}$, then the Leray spectral sequence ensures that $H^n(X,\ \mathbb{C}) = H^n(\mathbb{P}^n, \mathfrak{F})$ and $H^n(X^{'},\ \mathbb{C}) = H^n(\mathbb{P}^n, \mathfrak{F}^{'})$. Under the action of $\mathbb{Z}/r\mathbb{Z}$, we have the eigen subsheaf decomposition

$$\mathfrak{F} = \oplus_{j=0}^{r-1} \mathfrak{F}_{(j)}.$$

Similarly, we have

$$\mathfrak{F}^{'} = \oplus_{j=0}^{r^{'}-1} \mathfrak{F}^{'}_{(j)}$$

under the action of $\mathbb{Z}/r^{'}\mathbb{Z}$. It is not difficult to see the stalks $\mathfrak{F}_{(i),x} = \mathfrak{F}^{'}_{(1),x} = 0$, for each $x \in \cup_{j=1}^m H_j$. Let $U = \mathbb{P}^n - \cup_{j=1}^m H_j$, we can verify both $\mathfrak{F}_{(i)}|_U$ and $\mathfrak{F}^{'}_{(1)}|_U$ are rank one local systems on $U$. Moreover, we can show $\mathfrak{F}_{(i)}|_U \simeq \mathfrak{F}^{'}_{(1)}|_U$ by comparing the induced homomorphisms from $\pi_1(U)$ to $\mathbb{C}^*$. This completes the proof.

□

The following claim relates $\widetilde{Mon}_{(i)}$ and $Mon_{(i)}$.

**Claim 4.2.** For each integer $i$ with $1 \leq i \leq \llcorner \frac{r}{2} \lrcorner$, we have $\widetilde{Mon}^0_{(i)} = Mon^0_{(i)}$.

Indeed, we can verify without difficulty that, the monodromy transformation on $H^n(X,\ \mathbb{C})_{(i)}$ induced by any loop in the fiber of $\widetilde{\mathfrak{M}} \to \mathfrak{M}_{AR}$ is a $\zeta_r^j$-scalar multiplication, for some $j \in \mathbb{Z}$. Then the claim follows from the following commutative diagram

$$\begin{array}{ccc} \pi_1(\widetilde{\mathfrak{M}}, \tilde{s}) & \longrightarrow & \pi_1(\mathfrak{M}_{AR}, s) \\ & \searrow^{\tilde{\rho_i}} \quad {}^{\rho_i}\swarrow & \\ & PAut(H^n(X,\ \mathbb{C})_{(i)}) & \end{array}$$

4.5. **Proof of Theorem 1.1.** To simplify the notations, we use $\widetilde{Mon}$ to denote $\widetilde{Mon}_{(1)}$. Let $\widetilde{Mon}(\mathbb{C})$ be the group of complex points of $\widetilde{Mon}$. By Claim 4.2, it suffices to verify $\widetilde{Mon}(\mathbb{C}) = Aut(H^n(X,\ \mathbb{C})_{(1)}, Q)$.

Recall from subsection 4.1 that we have an embedding $\widetilde{\mathfrak{M}}_C \hookrightarrow \widetilde{\mathfrak{M}}$, the base point $\tilde{s} = (l_1, \cdots, l_m) \in \widetilde{\mathfrak{M}}_C$, and $C$ (resp. $X$) is the fiber over $\tilde{s}$ of the family $\widetilde{\mathcal{C}}$ (resp. $\widetilde{\mathcal{X}}$).

Obviously the permutation group $S_m$ acts naturally on $\widetilde{\mathfrak{M}}_C$ and $\widetilde{\mathfrak{M}}$. For each $\pi \in S_m$, we choose a path $\gamma_\pi$ in $\widetilde{\mathfrak{M}}_C$ from $\tilde{s}$ to $\pi(\tilde{s}) = (l_{\pi(1)}, \cdots, l_{\pi(m)})$. Let $\overline{\mathfrak{M}} = \widetilde{\mathfrak{M}}/S_m$ be the quotient space. Then we have group homomorphisms

$$\pi_1(\widetilde{\mathfrak{M}}_C, \tilde{s}) \to \pi_1(\widetilde{\mathfrak{M}}, \tilde{s}) \to \pi_1(\overline{\mathfrak{M}}, \bar{s}).$$

Here we use $\bar{s}$ to denote the image of $\tilde{s} \in \widetilde{\mathfrak{M}}$ in $\overline{\mathfrak{M}}$.

Since $\widetilde{\mathfrak{M}} \to \overline{\mathfrak{M}}$ is a Galois covering, $\pi_1(\widetilde{\mathfrak{M}}, \tilde{s})$ is a normal subgroup of $\pi_1(\overline{\mathfrak{M}}, \bar{s})$. We can verify without difficulty that $\pi_1(\overline{\mathfrak{M}}, \bar{s})$ is generated by $\pi_1(\widetilde{\mathfrak{M}}, \tilde{s})$ and $[\bar{\gamma}_\pi]$ $(\pi \in S_m)$, where $\bar{\gamma}_\pi$ is the loop in $\overline{\mathfrak{M}}$ represented by the image of $\gamma_\pi$. Moreover, since $\gamma_\pi$ is a path in $\widetilde{\mathfrak{M}}_C$, we know for each $\pi \in S_m$, the element $[\bar{\gamma}_\pi]$ normalizes the image of $\pi_1(\widetilde{\mathfrak{M}}_C, \tilde{s})$ in $\pi_1(\overline{\mathfrak{M}}, \bar{s})$.

Let $W = H^1(C,\ \mathbb{C})_{(1)} = H^1(C,\ \mathbb{C})$ and $V = H^n(X,\ \mathbb{C})_{(1)}$. Obviously the family $\widetilde{\mathcal{X}} \to \widetilde{\mathfrak{M}}$ descends to a family $\overline{\mathcal{X}} \to \overline{\mathfrak{M}}$, and we have the monodromy representation

$$\bar{\rho} : \pi_1(\overline{\mathfrak{M}}, \bar{s}) \to Aut(V).$$

Let $\overline{Mon}$ be the smallest algebraic subgroup of $Aut(V)$ containing the image of $\pi_1(\overline{\mathfrak{M}}, \bar{s})$. By Proposition 4.1 and the commutative diagram (11), we have $V \simeq \wedge^n W$ and $Sp(W)$ is the smallest algebraic subgroup of $Aut(V)$ containing the image of $\pi_1(\widetilde{\mathfrak{M}}_C, \tilde{s})$. Then we get a sequence of group embeddings

$$Sp(W) \hookrightarrow \widetilde{Mon}(\mathbb{C}) \hookrightarrow \overline{Mon} \hookrightarrow Aut(V) = Aut(\wedge^n W)$$

satisfying: $\overline{Mon}$ is generated by $\widetilde{Mon}(\mathbb{C})$ and $\bar{\rho}([\bar{\gamma}_\pi])$ $(\pi \in S_m)$. Moreover, $\bar{\rho}([\bar{\gamma}_\pi])$ normalizes $Sp(W)$, for each $\pi \in S_m$.

We first claim: any $\widetilde{Mon}(\mathbb{C})$-invariant linear subspace $U$ of $V$ is also $\overline{Mon}$-invariant.

In fact, since $\overline{Mon}$ is generated by $\widetilde{Mon}(\mathbb{C})$ and $\bar{\rho}([\bar{\gamma}_\pi])$ $(\pi \in S_m)$, we only need to verify $U$ is invariant under the action of $\bar{\rho}([\bar{\gamma}_\pi])$, for each $\pi \in S_m$. This follows from the elementary Lemma 4.8 below.

By Proposition 4.2, the set $R_{\mathbb{C}} := \{\gamma \cdot e_{k_1k_2\cdots k_{n+1},(1)} | \gamma \in \widetilde{Mon}(\mathbb{C}),\ 1 \leq k_1 < k_2 < \cdots < k_{n+1} \leq m\}$ linearly spans $V$ and by definition, $R_{\mathbb{C}}$ consists of a finite union of $\widetilde{Mon}(\mathbb{C})$-orbits.

We next claim: $\forall\ \alpha \neq 0 \in R_{\mathbb{C}}$, the orbit $\overline{Mon} \cdot \alpha$ linearly spans $V$.

In fact, given $1 \leq k_1 < \cdots < k_{n+1} \leq m$ and $1 \leq l_1 < \cdots < l_{n+1} \leq m$, we can choose a permutation $\pi \in S_m$ such that $\pi(k_i) = l_i$, for each $1 \leq i \leq n+1$. Recall $\gamma_{k_1\cdots k_{n+1}}$ is a meridian of $\Delta_{k_1\cdots k_{n+1}}$ with base point $\tilde{s}$, so $\pi(\gamma_{k_1\cdots k_{n+1}})$ is a meridian of $\Delta_{l_1\cdots l_{n+1}}$ with base point $\pi(\tilde{s})$. Take any path $\gamma_1$ in $\widetilde{\mathfrak{M}}$ from $\tilde{s}$ to $\pi(\tilde{s})$, then $\gamma_1 \circ \pi(\gamma_{k_1\cdots k_{n+1}}) \circ \gamma_1^{-1}$ is a meridian of $\Delta_{l_1\cdots l_{n+1}}$ with base point $\tilde{s}$. Note the image of $\gamma_1$ becomes a loop with base point $\bar{s}$ under the covering map $\widetilde{\mathfrak{M}} \to \overline{\mathfrak{M}}$. We find the monodromy actions $\bar{\rho}(\gamma_{k_1\cdots k_{n+1}})$ and $\bar{\rho}(\gamma_{l_1\cdots l_{n+1}})$ are conjugate in $\overline{Mon}$. Suppose $g \in \overline{Mon}$ satisfying $g\bar{\rho}(\gamma_{k_1\cdots k_{n+1}})g^{-1} = \bar{\rho}(\gamma_{l_1\cdots l_{n+1}})$. By the Picard-Lefschetz type formula (12), this in turn implies that there exists a nonzero $\lambda \in \mathbb{C}$, such that $ge_{k_1\cdots k_{n+1},(1)} = \lambda\ e_{l_1\cdots l_{n+1},(1)}$. Since $R_{\mathbb{C}}$ linearly spans $V$, our claim follows.

We apply Proposition 4.3 and Proposition 4.4 to $\widetilde{Mon}(\mathbb{C})$. Clearly $\widetilde{Mon}(\mathbb{C})$ is not a finite group. If $\widetilde{Mon}(\mathbb{C}) \neq Aut(V, Q)$, there would exist a nontrivial $\widetilde{Mon}(\mathbb{C})$-invariant subspace $U$ of $V$ such that $R_{\mathbb{C}} \subset U \cup U^{\perp}$. Then both $U$ and $U^{\perp}$ are also $\overline{Mon}$-invariant. Now that $\forall\ \alpha \neq 0 \in R_{\mathbb{C}}$, the orbit $\overline{Mon} \cdot \alpha$ linearly spans $V$, it follows that $U = V$ or $U^{\perp} = V$. This contradicts with the nontrivial assumption of $U$. So we get $\widetilde{Mon}(\mathbb{C}) = Aut(V, Q)$. This completes the proof of Theorem 1.1.

**Lemma 4.8.** Suppose $W$ is a finite dimensional complex vector space and $Q$ is a non-degenerate alternating bilinear form on $W$. Denote the symplectic group $Sp(W, Q)$ by $Sp(W)$. Let $V = \wedge^n W$, with $n \leq dimW - 1$. We have a natural embedding $Sp(W) \hookrightarrow Aut(V) = Aut(\wedge^n W)$ through the $Sp(W)$-action on $V$. Suppose $g \in Aut(V)$ such that $g$ normalizes $Sp(W)$. Then any $Sp(W)$-invariant linear subspace $U$ of $V$ is also $g$-invariant.

*Proof.* Since $g$ normalizes $Sp(W)$, for any $Sp(W)$-invariant linear subspace $U$ of $V$, the subspace $gU \subset V$ is also $Sp(W)$-invariant. It is well known that any automorphism of the complex linear algebraic group $Sp(W)$ is an inner automorphism (cf. [12]). So there exists $h \in Sp(W)$ such that $\forall\ x \in Sp(W)$, we have $gh^{-1}xhg^{-1} = x$. Then $U \xrightarrow{gh^{-1}} gU$ is an isomorphism of $Sp(W)$-modules. On the other hand, let $V = \oplus_{i=1}^{s} V_i$ be the irreducible decomposition as a $Sp(W)$-module, then we know for any $1 \leq i < j \leq s$, the $Sp(W)$-modules $V_i$ and $V_j$ are not isomorphic (cf. [9], Theorem 17.5). From this one can deduce $gU = U$ from the isomorphism $gU \simeq U$.

□

4.6. **Proof of Theorem 2.1.** By Claim 4.1, Theorem 2.1, (1) can be reduced to Theorem 1.1, which has been proven. In order to prove Theorem 2.1, (2), we assume $r > 2$, $1 \leq i < \frac{r}{2}$ and $mi \geq 2r$. By Claim 4.2, it suffices to prove $\widetilde{Mon}^0_{(i)} = SU(H^n(X,\ \mathbb{C})_{(i)}, H)$. Let $W = H^1(C,\ \mathbb{C})_{(i)}$ and $V = H^n(X,\ \mathbb{C})_{(i)} \simeq \wedge^n W$. By Proposition 4.1 and the commutative diagram (11), we have the embeddings

$$SU(W) \hookrightarrow \widetilde{Mon}_{(i)} \hookrightarrow Aut(V) = Aut(\wedge^n W).$$

By Proposition 4.2, the set $R$ linearly spans $V$ and it consists of a finite union of $\widetilde{Mon}_{(i)}$-orbits.

We divide the proof into two cases.

Case 1. $r \nmid (n+1)i$. In this case, Proposition 3.5 implies that for each $e \in R$, the intersection pairing $H(e,e) = (\zeta_r^{(n+1)i} - 1)c^{-1} \neq 0$. Since the representation of $SU(W)$ on $V = \wedge^n W$ is irreducible, we deduce $\widetilde{Mon}_{(i)}$ acts irreducibly on $V$, and hence there does not exist any nontrivial $\widetilde{Mon}_{(i)}$-invariant subspace of $V$. So we can apply Proposition 4.5 to obtain $P\widetilde{Mon}_{(i)} = PU(V,H)$, and this implies $\widetilde{Mon}^0_{(i)} = SU(V,H)$.

Case 2. $r|(n+1)i$. In this case, Proposition 3.5 implies that for each $e \in R$, the intersection pairing $H(e,e) = 0$. Taking a nonzero $e \in R$, the map

$$\begin{aligned} \Phi :& V \to V \\ & \alpha \mapsto \alpha + cH(\alpha, e)e \end{aligned}$$

belongs to $\widetilde{Mon}_{(i)}$ by the Picard-Lefschetz type formula (12). Since $H(e,e) = 0$ and $H$ is non-degenerate on $V$, we see for any positive integer $N$, $\Phi^N$ is a nontrivial pseudo-reflection on $V$. Taking a sufficiently large integer $N$ such that $\Phi^N \in \widetilde{Mon}^0_{(i)}$, we obtain an element in $\widetilde{Mon}^0_{(i)}$ which is a nontrivial pseudo-reflection on $V$. By a result of Deligne (cf. Corollary 4.2.9 in [4]), the group $\widetilde{Mon}_{(i)}$ is semi-simple. Now an application of Proposition 4.6 gives $\widetilde{Mon}^0_{(i)}(\mathbb{C}) = SL(V) = SU(V,H)(\mathbb{C})$. So $\widetilde{Mon}^0_{(i)} = SU(V,H)$. This completes the proof of Theorem 2.1.

**Remark 4.1.** In fact, our arguments give an independent proof of Proposition 4.1, based on the irreducibility of the monodromy representation (see e.g. [14], Proposition 5.1).

## 5. Funding

This work was supported by the National Natural Science Foundation of China [11301496, 11471298].

## Acknowledgements

I would like to thank Prof. Mao Sheng for generous encouragement during this work. Particular thanks go to Prof. D. T. Lê, who carefully read the manuscript and made numerous helpful suggestions. Finally, I would like to thank the anonymous referee for making several helpful suggestions for improving the exposition.

School of Mathematical Sciences, University of Science and Technology of China, No. 96 Jinzhai Road, Hefei, Anhui 230026, P. R. China
*Email address*: xujx02@ustc.edu.cn

# CORRIGENDUM TO “ZARISKI DENSITY OF MONODROMY GROUPS VIA A PICARD–LEFSCHETZ TYPE FORMULA”

JINXING XU

Abstract. We correct the proof of Proposition 4.2 in [Xu18], which used a false assertion about the representation of a symplectic group on an exterior power. We show that the relevant vanishing cohomology classes, together with their orbits under the actual image of the monodromy representation, span each nontrivial eigenspace of the middle cohomology. The argument does not require the condition $mi \geq 2r$. The statements of the main theorems are unchanged.

In this corrigendum, we correct an error in the proof of Proposition 4.2 in [Xu18]. All references to numbered statements and sections of [Xu18] refer to the published version. The original proof incorrectly asserted that the representation

$$\mathrm{Sp}(H^1(C,\mathbb{R})_{(1)}) \curvearrowright \bigwedge^n H^1(C,\mathbb{R})_{(1)}$$

has no nontrivial invariant line. Indeed, for a symplectic vector space $(W,\omega)$ and $n = 2k \leq \dim W$, the inverse symplectic form defines a bivector $\omega^{-1} \in \bigwedge^2 W$, and $(\omega^{-1})^{\wedge k}$ is a nonzero $\mathrm{Sp}(W,\omega)$-invariant vector in $\bigwedge^{2k} W$. The proof below avoids the incorrect assertion altogether. It strengthens Proposition 4.2 by removing the condition $mi \geq 2r$ and by using orbits under the monodromy group itself, rather than its algebraic closure. The statements of the main theorems are unchanged.

We briefly recall the setup introduced in [Xu18, §2.1 and §4.1]. Let $m, n, r$ be positive integers satisfying $r \mid m$ and $m \geq n+3$. Put $d := m/r$ and $V := \mathbb{C}^{n+1}$. Regard the parameters $v_i$ as nonzero elements of $V^\vee$. Thus a point of the parameter space is $\underline{v} = (v_1, \ldots, v_m) \in (V^\vee \backslash \{0\})^m$. Choose homogeneous coordinates $[x_0 : \cdots : x_n]$ on $\mathbb{P}^n$ and write $v_i = \ell_i := a_{i0}x_0 + \cdots + a_{in}x_n$.

Let $\mathcal{I}$ be the set of tuples $I = (k_1, \ldots, k_{n+1})$ satisfying $1 \leq k_1 < k_2 < \cdots < k_{n+1} \leq m$. For $I \in \mathcal{I}$, let $\Delta_I$ be the locus where the linear forms $v_{k_1}, \ldots, v_{k_{n+1}}$ are linearly dependent in $V^\vee$. Each $\Delta_I$ is an irreducible hypersurface in $(V^\vee \setminus \{0\})^m$. Define

$$\widetilde{\mathfrak{M}} := (V^\vee \setminus \{0\})^m \setminus \bigcup_{I \in \mathcal{I}} \Delta_I.$$

For $\underline{v} \in \widetilde{\mathfrak{M}}$, the hyperplane $H_{v_i} := \mathbb{P}(\ker v_i) \subset \mathbb{P}^n$ is defined by the linear form $\ell_i$. Let $\pi\colon \mathrm{Tot}(\mathcal{O}_{\mathbb{P}^n}(d)) \longrightarrow \mathbb{P}^n$ be the total space of $\mathcal{O}_{\mathbb{P}^n}(d)$.

The degree-$r$ cyclic cover branched along $\sum_{i=1}^m H_{v_i}$ is the hypersurface $X_{\underline{v}}$ defined by

$$y^r - \pi^*(\ell_1 \cdots \ell_m) = 0$$

in $\mathrm{Tot}(\mathcal{O}_{\mathbb{P}^n}(d))$, where $y$ is the tautological section of $\pi^*\mathcal{O}_{\mathbb{P}^n}(d)$. More precisely,

$$y^r,\ \pi^*(\ell_1 \cdots \ell_m) \in H^0\big(\mathrm{Tot}(\mathcal{O}_{\mathbb{P}^n}(d)), \pi^*\mathcal{O}_{\mathbb{P}^n}(m)\big),$$

because $rd = m$.

Fix $\zeta_r \coloneqq \exp(2\pi \mathrm{i}/r)$. The automorphism

$$([x_0 : \cdots : x_n], y) \longmapsto ([x_0 : \cdots : x_n], \zeta_r y)$$

of $\mathrm{Tot}(\mathcal{O}_{\mathbb{P}^n}(d))$ induces an automorphism $\sigma$ of $X_{\underline{v}}$.

As in [Xu18, §4.1], these covers form a topologically locally trivial family $\tilde{f}\colon \widetilde{\mathcal{X}} \longrightarrow \widetilde{\mathfrak{M}}$ whose fiber over $\underline{v}$ is $X_{\underline{v}}$. The automorphism $\sigma$ extends naturally to this family.

Fix a base point $\tilde{s} \in \widetilde{\mathfrak{M}}$ and set $X \coloneqq \tilde{f}^{-1}(\tilde{s})$. For $1 \le i \le r-1$, put

$$V_i \coloneqq H^n(X,\mathbb{C})_{(i)} \coloneqq \left\{ \alpha \in H^n(X,\mathbb{C}) \mid \sigma^*\alpha = \zeta_r^i \alpha \right\}.$$

Let $Q$ denote the cup-product pairing on $H^n(X,\mathbb{C})$. We use the non-degenerate Hermitian form

$$H(\alpha,\beta) \coloneqq \mathrm{i}^n Q(\alpha,\overline{\beta}),$$

which is non-degenerate on $V_i$ and invariant under monodromy; these properties are established in [Xu18, §2.2–2.3]. Here $H$ is linear in its first variable and conjugate-linear in its second variable. Let $U(V_i,H)$ denote the real algebraic group of complex linear automorphisms of $V_i$ preserving $H$. Consider the cohomological monodromy representation $\tilde{\rho}_i\colon \pi_1(\widetilde{\mathfrak{M}},\tilde{s}) \longrightarrow U(V_i,H) \subset \mathrm{GL}_{\mathbb{C}}(V_i)$, and write $\Gamma_i \coloneqq \mathrm{Im}\,\tilde{\rho}_i$.

By the Picard–Lefschetz analysis in [Xu18, §3.3 and Proposition 3.5], for each $I \in \mathcal{I}$ we may choose a meridian $\gamma_I$ of $\Delta_I$ and a corresponding vanishing cohomology class $e_{I,(i)} \in V_i$ such that

$$\tilde{\rho}_i(\gamma_I)(\alpha) = \alpha + c_{I,i} H(\alpha, e_{I,(i)}) e_{I,(i)}, \qquad \alpha \in V_i. \tag{0.1}$$

Here $c_{I,i} \in \mathbb{C}^\times$ is a nonzero scalar. The class $e_{I,(i)}$ is the cohomology class Poincaré dual to the corresponding vanishing cycle. No nonvanishing assumption on these classes is needed in the argument below.

Define

$$R_i^{\mathrm{mon}} \coloneqq \left\{ g(e_{I,(i)}) \mid g \in \Gamma_i,\ I \in \mathcal{I} \right\} \subset V_i.$$

**Proposition 0.1** (Refinement of Proposition 4.2 of [Xu18])**.** *For every $1 \le i \le r-1$, the $\mathbb{C}$-linear span of $R_i^{\mathrm{mon}}$ is $V_i$. In particular, no condition $mi \ge 2r$ is required.*

*Proof.* Consider the set

$$S \coloneqq \left\{ \gamma \circ \gamma_I \circ \gamma^{-1} \mid \gamma \in \pi_1(\widetilde{\mathfrak{M}},\tilde{s}),\ I \in \mathcal{I} \right\}.$$

Since $n \geq 1$, the ambient space $(V^\vee \setminus \{0\})^m$ is simply connected. The meridians around the irreducible components $\Delta_I$ normally generate $\pi_1(\widetilde{\mathfrak{M}}, \tilde{s})$. Indeed, a loop bounds a disk in the ambient space; after a general-position perturbation, the disk avoids the singular locus of the divisor and meets its smooth locus transversely in finitely many points. Removing small disks around these points expresses the loop as a product of conjugates of meridians and their inverses. Because $S$ contains all their conjugates, $S$ generates $\pi_1(\widetilde{\mathfrak{M}}, \tilde{s})$.

Because $H$ is invariant under monodromy, (0.1) gives, for $\gamma \in \pi_1(\widetilde{\mathfrak{M}}, \tilde{s})$ and $I \in \mathcal{I}$,

$$\tilde{\rho}_i(\gamma \circ \gamma_I \circ \gamma^{-1})(\alpha) = \alpha + c_{I,i} H(\alpha, \tilde{\rho}_i(\gamma) e_{I,(i)}) \tilde{\rho}_i(\gamma) e_{I,(i)}. \tag{0.2}$$

Suppose, for contradiction, that $R_i^{\mathrm{mon}}$ does not span $V_i$. Since $H$ is nondegenerate, there exists a nonzero element $\alpha_0 \in V_i$ such that $H(\alpha_0, e) = 0$ for all $e \in R_i^{\mathrm{mon}}$. For every element of $S$, the second term in (0.2) therefore vanishes. Since $S$ generates the fundamental group, $\alpha_0$ is fixed by $\Gamma_i$.

Write $\tilde{s} = (v_1, \ldots, v_m)$ and consider the loop

$$\gamma_0 \colon [0,1] \longrightarrow \widetilde{\mathfrak{M}}, \qquad t \longmapsto (\exp(-2\pi \mathrm{i} t) v_1, v_2, \ldots, v_m).$$

Scaling $v_1$ by a nonzero scalar does not change any linear-dependence relation, so $\gamma_0$ is a loop in $\widetilde{\mathfrak{M}}$. We now compute its cohomological monodromy action: $\tilde{\rho}_i(\gamma_0)(\alpha) = \zeta_r^i \alpha$ for all $\alpha \in V_i$. Thus $\tilde{\rho}_i(\gamma_0)(\alpha_0) = \zeta_r^i \alpha_0$, contradicting the fact that $\alpha_0$ is fixed by $\Gamma_i$, because $\zeta_r^i \neq 1$.

To verify the displayed monodromy action, consider the pullback family $\gamma_0^* \widetilde{\mathcal{X}} \to [0,1]$. Using coordinates $([x_0 : \cdots : x_n], y)$ on $\mathrm{Tot}(\mathcal{O}_{\mathbb{P}^n}(d))$, define

$$\begin{aligned} \varphi \colon [0,1] \times X &\xrightarrow{\sim} \gamma_0^* \widetilde{\mathcal{X}}, \\ (t, ([x_0 : \cdots : x_n], y)) &\longmapsto \left( t, [x_0 : \cdots : x_n], \exp\left( -\frac{2\pi \mathrm{i} t}{r} \right) y \right). \end{aligned}$$

At $t = 1$, the endpoint map is $\varphi_1 = \sigma^{-1}$. Since $\tilde{\rho}_i$ denotes the monodromy of the cohomology local system, its action is

$$\tilde{\rho}_i(\gamma_0) = (\varphi_1^{-1})^* = \sigma^*.$$

By the definition of $V_i$, this is multiplication by $\zeta_r^i$, as claimed. □

Let $R_i^{\mathrm{alg}}$ denote the set $R$ defined in [Xu18, §4.1], using the algebraic monodromy group. Since $\Gamma_i$ is contained in that group, we have $R_i^{\mathrm{mon}} \subseteq R_i^{\mathrm{alg}}$. Consequently, Proposition 0.1 implies the spanning assertion of [Xu18, Proposition 4.2], without the hypothesis $mi \geq 2r$. The original sets of algebraic-monodromy orbits can therefore be retained in the subsequent arguments.

*Remark* 0.2 (Fermat subfamilies and vanishing cycles of hypersurfaces). The same argument gives a direct proof of the classical spanning property for vanishing cycles of smooth hypersurfaces in projective space; see, for example,

[Bea86, §1]. Let

$$U \subset \mathbb{P}H^0(\mathbb{P}^{n+1}, \mathcal{O}(d))$$

be the parameter space of smooth degree-$d$ hypersurfaces, where in this remark $n \geq 1$ and $d \geq 2$ are arbitrary. Take the Fermat hypersurface

$$F = \{x_0^d + x_1^d + \cdots + x_{n+1}^d = 0\}$$

as the base fiber. Put $P = H^n(F, \mathbb{Q})_{\mathrm{prim}}$, and let $\Gamma$ be the image of the monodromy representation of $\pi_1(U, [F])$ on $P$. Let $E \subseteq P$ be the $\mathbb{Q}$-linear span of the nodal vanishing cohomology classes of the universal family, transported to $F$ along all possible paths.

The discriminant is an irreducible hypersurface. Its general point represents a hypersurface with a single ordinary double point. Its meridians normally generate $\pi_1(U, [F])$, since the ambient parameter space is simply connected. The Picard–Lefschetz formula therefore gives

$$E^{\perp} = P^{\Gamma},$$

where the orthogonal complement is taken in $P$ for the cup-product pairing; compare [Del74, Proposition 5.3].

To show that $P^{\Gamma} = 0$, consider the Fermat subfamily

$$F_a = \{a x_0^d + x_1^d + \cdots + x_{n+1}^d = 0\}, \qquad a \in \mathbb{C}^{\times}.$$

Write $\zeta_d = \exp(2\pi \mathrm{i}/d)$ and set

$$\sigma([x_0 : x_1 : \cdots : x_{n+1}]) = [\zeta_d x_0 : x_1 : \cdots : x_{n+1}].$$

Along the loop $a(t) = \exp(2\pi \mathrm{i} t)$, a trivialization is given by

$$\varphi_t : F \longrightarrow F_{a(t)},$$
$$[x_0 : x_1 : \cdots : x_{n+1}] \longmapsto [\exp(-2\pi \mathrm{i} t/d) x_0 : x_1 : \cdots : x_{n+1}].$$

Since $\varphi_1 = \sigma^{-1}$, its cohomological monodromy is $(\varphi_1^{-1})^* = \sigma^*$. In particular, $\sigma^*|_P \in \Gamma$.

The projection

$$p : F \longrightarrow \mathbb{P}^n, \qquad [x_0 : x_1 : \cdots : x_{n+1}] \longmapsto [x_1 : \cdots : x_{n+1}]$$

identifies $\mathbb{P}^n$ with $F/\langle \sigma \rangle$. Thus

$$H^n(F, \mathbb{Q})^{\sigma^*} = p^* H^n(\mathbb{P}^n, \mathbb{Q}).$$

Writing $h = c_1(\mathcal{O}_F(1))$, we have $p^* \mathcal{O}_{\mathbb{P}^n}(1) = \mathcal{O}_F(1)$, and hence

$$H^n(F, \mathbb{Q})^{\sigma^*} = \begin{cases} 0, & n \text{ odd}, \\ \mathbb{Q} h^{n/2}, & n \text{ even}. \end{cases}$$

Thus the invariant part has trivial intersection with the primitive cohomology, and consequently $P^{\sigma^*} = 0$. Therefore

$$E^{\perp} = P^{\Gamma} \subseteq P^{\sigma^*} = 0.$$

The nondegeneracy of the primitive cup-product pairing gives $E = P$.

## Acknowledgments

The author thanks Ben Castor, Haohua Deng, Matt Kerr, and Gregory Pearlstein for their paper [CDKP23], which brought to his attention the error in [Xu18] discussed above.

School of Mathematical Sciences, University of Science and Technology of China, No. 96 Jinzhai Road, Hefei, Anhui 230026, P. R. China

*Email address*: xujx02@ustc.edu.cn